%

\catcode`\@=11
\font\tensmc=cmcsc10      
\def\smc{\tensmc}

\def\hcorrection#1{\advance\hoffset by #1 }
\def\vcorrection#1{\advance\voffset by #1 }
\def\wlog#1{}
\newif\iftitle@
\outer\def\title{\title@true\vglue 24\p@ plus 12\p@ minus 12\p@
   \bgroup\let\\=\cr\tabskip\centering
   \halign to \hsize\bgroup\tenbf\hfill\ignorespaces##\unskip\hfill\cr}
\def\endtitle{\cr\egroup\egroup\vglue 18\p@ plus 12\p@ minus 6\p@}
\outer\def\author{\iftitle@\vglue -18\p@ plus -12\p@ minus -6\p@\fi\vglue
    12\p@ plus 6\p@ minus 3\p@\bgroup\let\\=\cr\tabskip\centering
    \halign to \hsize\bgroup\smc\hfill\ignorespaces##\unskip\hfill\cr}
\def\endauthor{\cr\egroup\egroup\vglue 18\p@ plus 12\p@ minus 6\p@}
\outer\def\heading{\bigbreak\bgroup\let\\=\cr\tabskip\centering
    \halign to \hsize\bgroup\smc\hfill\ignorespaces##\unskip\hfill\cr}
\def\endheading{\cr\egroup\egroup\nobreak\medskip}

\outer\def\endproclaim{\par\ifdim\lastskip<\medskipamount\removelastskip
  \penalty 55 \fi\medskip\rm}
\outer\def\demo#1{\par\ifdim\lastskip<\smallskipamount\removelastskip
    \smallskip\fi\noindent{\smc\ignorespaces#1\unskip:\enspace}\rm
      \ignorespaces}

\newcount\footmarkcount@
\footmarkcount@=1
\def\makefootnote@#1#2{\insert\footins{\interlinepenalty=100
  \splittopskip=\ht\strutbox \splitmaxdepth=\dp\strutbox 
  \floatingpenalty=\@MM
  \leftskip=\z@\rightskip=\z@\spaceskip=\z@\xspaceskip=\z@
  \noindent{#1}\footstrut\rm\ignorespaces #2\strut}}
\def\footnote{\let\@sf=\empty\ifhmode\edef\@sf{\spacefactor
   =\the\spacefactor}\/\fi\futurelet\next\footnote@}
\def\footnote@{\ifx"\next\let\next\footnote@@\else
    \let\next\footnote@@@\fi\next}
\def\footnote@@"#1"#2{#1\@sf\relax\makefootnote@{#1}{#2}}
\def\footnote@@@#1{$^{\number\footmarkcount@}$\makefootnote@
   {$^{\number\footmarkcount@}$}{#1}\global\advance\footmarkcount@ by 1 }

\hyphenation{man-u-script man-u-scripts ap-pen-dix ap-pen-di-ces}
\hyphenation{data-base data-bases}
\ifx\amstexloaded@\relax\catcode`\@=13 
  \endinput\else\let\amstexloaded@=\relax\fi
\newlinechar=`\^^J
\def\eat@#1{}
\def\Space@.{\futurelet\Space@\relax}
\Space@. %
\newhelp\athelp@
{Only certain combinations beginning with @ make sense to me.^^J
Perhaps you wanted \string\@\space for a printed @?^^J
I've ignored the character or group after @.}
\def\futureletnextat@{\futurelet\next\at@}
{\catcode`\@=\active
\lccode`\Z=`\@ \lowercase
{\gdef@{\expandafter\csname futureletnextatZ\endcsname}
\expandafter\gdef\csname atZ\endcsname
   {\ifcat\noexpand\next a\def\next{\csname atZZ\endcsname}\else
   \ifcat\noexpand\next0\def\next{\csname atZZ\endcsname}\else
    \def\next{\csname atZZZ\endcsname}\fi\fi\next}
\expandafter\gdef\csname atZZ\endcsname#1{\expandafter
   \ifx\csname #1Zat\endcsname\relax\def\next
     {\errhelp\expandafter=\csname athelpZ\endcsname
      \errmessage{Invalid use of \string@}}\else
       \def\next{\csname #1Zat\endcsname}\fi\next}
\expandafter\gdef\csname atZZZ\endcsname#1{\errhelp
    \expandafter=\csname athelpZ\endcsname
      \errmessage{Invalid use of \string@}}}}
\def\atdef@#1{\expandafter\def\csname #1@at\endcsname}
\newhelp\defahelp@{If you typed \string\define\space cs instead of
\string\define\string\cs\space^^J
I've substituted an inaccessible control sequence so that your^^J
definition will be completed without mixing me up too badly.^^J
If you typed \string\define{\string\cs} the inaccessible control sequence^^J
was defined to be \string\cs, and the rest of your^^J
definition appears as input.}
\newhelp\defbhelp@{I've ignored your definition, because it might^^J
conflict with other uses that are important to me.}
\def\define{\futurelet\next\define@}
\def\define@{\ifcat\noexpand\next\relax
  \def\next{\define@@}%
  \else\errhelp=\defahelp@
  \errmessage{\string\define\space must be followed by a control 
     sequence}\def\next{\def\garbage@}\fi\next}
\def\undefined@{}
\def\preloaded@{}    
\def\define@@#1{\ifx#1\relax\errhelp=\defbhelp@
   \errmessage{\string#1\space is already defined}\def\next{\def\garbage@}%
   \else\expandafter\ifx\csname\expandafter\eat@\string
         #1@\endcsname\undefined@\errhelp=\defbhelp@
   \errmessage{\string#1\space can't be defined}\def\next{\def\garbage@}%
   \else\expandafter\ifx\csname\expandafter\eat@\string#1\endcsname\relax
     \def\next{\def#1}\else\errhelp=\defbhelp@
     \errmessage{\string#1\space is already defined}\def\next{\def\garbage@}%
      \fi\fi\fi\next}
\def\famzero{\fam\z@}

\def\lim{\mathop{\famzero lim}}

\def\min{\mathop{\famzero min}}

\def\sup{\mathop{\famzero sup}}

\def\textfont@#1#2{\def#1{\relax\ifmmode
    \errmessage{Use \string#1\space only in text}\else#2\fi}}
\textfont@\rm\tenrm
\textfont@\it\tenit
\textfont@\sl\tensl
\textfont@\bf\tenbf
\textfont@\smc\tensmc
\let\ic@=\/
\def\/{\unskip\ic@}
\def\textfonti{\the\textfont1 }
\def\t#1#2{{\edef\next{\the\font}\textfonti\accent"7F \next#1#2}}
\let\B=\=
\let\D=\.
\def~{\unskip\nobreak\ \ignorespaces}
{\catcode`\@=\active
\gdef\@{\char'100 }}
\atdef@-{\leavevmode\futurelet\next\athyph@}
\def\athyph@{\ifx\next-\let\next=\athyph@@
  \else\let\next=\athyph@@@\fi\next}
\def\athyph@@@{\hbox{-}}
\def\athyph@@#1{\futurelet\next\athyph@@@@}
\def\athyph@@@@{\if\next-\def\next##1{\hbox{---}}\else
    \def\next{\hbox{--}}\fi\next}
\def\.{.\spacefactor=\@m}
\atdef@.{\null.}
\atdef@,{\null,}
\atdef@;{\null;}
\atdef@:{\null:}
\atdef@?{\null?}
\atdef@!{\null!}   
\def\srdr@{\thinspace}                     
\def\drsr@{\kern.02778em}
\def\sldl@{\kern.02778em}
\def\dlsl@{\thinspace}
\atdef@"{\unskip\futurelet\next\atqq@}
\def\atqq@{\ifx\next\Space@\def\next. {\atqq@@}\else
         \def\next.{\atqq@@}\fi\next.}
\def\atqq@@{\futurelet\next\atqq@@@}
\def\atqq@@@{\ifx\next`\def\next`{\atqql@}\else\def\next'{\atqqr@}\fi\next}
\def\atqql@{\futurelet\next\atqql@@}
\def\atqql@@{\ifx\next`\def\next`{\sldl@``}\else\def\next{\dlsl@`}\fi\next}
\def\atqqr@{\futurelet\next\atqqr@@}
\def\atqqr@@{\ifx\next'\def\next'{\srdr@''}\else\def\next{\drsr@'}\fi\next}

\def\textfontii{\the\textfont2 }
\def\{{\relax\ifmmode\lbrace\else
    {\textfontii f}\spacefactor=\@m\fi}
\def\}{\relax\ifmmode\rbrace\else
    \let\@sf=\empty\ifhmode\edef\@sf{\spacefactor=\the\spacefactor}\fi
      {\textfontii g}\@sf\relax\fi}   
\def\nonhmodeerr@#1{\errmessage
     {\string#1\space allowed only within text}}
\def\linebreak{\relax\ifhmode\unskip\break\else
    \nonhmodeerr@\linebreak\fi}
\def\allowlinebreak{\relax
   \ifhmode\allowbreak\else\nonhmodeerr@\allowlinebreak\fi}
\newskip\saveskip@
\def\nolinebreak{\relax\ifhmode\saveskip@=\lastskip\unskip
  \nobreak\ifdim\saveskip@>\z@\hskip\saveskip@\fi
   \else\nonhmodeerr@\nolinebreak\fi}
\def\newline{\relax\ifhmode\null\hfil\break
    \else\nonhmodeerr@\newline\fi}
\def\nonmathaerr@#1{\errmessage
     {\string#1\space is not allowed in display math mode}}
\def\nonmathberr@#1{\errmessage{\string#1\space is allowed only in math mode}}
\def\mathbreak{\relax\ifmmode\ifinner\break\else
   \nonmathaerr@\mathbreak\fi\else\nonmathberr@\mathbreak\fi}
\def\nomathbreak{\relax\ifmmode\ifinner\nobreak\else
    \nonmathaerr@\nomathbreak\fi\else\nonmathberr@\nomathbreak\fi}
\def\allowmathbreak{\relax\ifmmode\ifinner\allowbreak\else
     \nonmathaerr@\allowmathbreak\fi\else\nonmathberr@\allowmathbreak\fi}
\def\pagebreak{\relax\ifmmode
   \ifinner\errmessage{\string\pagebreak\space
     not allowed in non-display math mode}\else\postdisplaypenalty-\@M\fi
   \else\ifvmode\penalty-\@M\else\edef\spacefactor@
       {\spacefactor=\the\spacefactor}\vadjust{\penalty-\@M}\spacefactor@
        \relax\fi\fi}
\def\nopagebreak{\relax\ifmmode
     \ifinner\errmessage{\string\nopagebreak\space
    not allowed in non-display math mode}\else\postdisplaypenalty\@M\fi
    \else\ifvmode\nobreak\else\edef\spacefactor@
        {\spacefactor=\the\spacefactor}\vadjust{\penalty\@M}\spacefactor@
         \relax\fi\fi}
\def\newpage{\relax\ifvmode\vfill\penalty-\@M\else\nonvmodeerr@\newpage\fi}
\def\nonvmodeerr@#1{\errmessage
    {\string#1\space is allowed only between paragraphs}}
\def\smallpagebreak{\relax\ifvmode\smallbreak
      \else\nonvmodeerr@\smallpagebreak\fi}
\def\medpagebreak{\relax\ifvmode\medbreak
       \else\nonvmodeerr@\medpagebreak\fi}
\def\bigpagebreak{\relax\ifvmode\bigbreak
      \else\nonvmodeerr@\bigpagebreak\fi}
\newdimen\captionwidth@
\captionwidth@=\hsize
\advance\captionwidth@ by -1.5in
\def\caption#1{}
\def\topspace#1{\gdef\thespace@{#1}\ifvmode\def\next
    {\futurelet\next\topspace@}\else\def\next{\nonvmodeerr@\topspace}\fi\next}
\def\topspace@{\ifx\next\Space@\def\next. {\futurelet\next\topspace@@}\else
     \def\next.{\futurelet\next\topspace@@}\fi\next.}
\def\topspace@@{\ifx\next\caption\let\next\topspace@@@\else
    \let\next\topspace@@@@\fi\next}
 \def\topspace@@@@{\topinsert\vbox to 
       \thespace@{}\endinsert}
\def\topspace@@@\caption#1{\topinsert\vbox to
    \thespace@{}\nobreak
      \smallskip
    \setbox\z@=\hbox{\noindent\ignorespaces#1\unskip}%
   \ifdim\wd\z@>\captionwidth@
   \centerline{\vbox{\hsize=\captionwidth@\noindent\ignorespaces#1\unskip}}%
   \else\centerline{\box\z@}\fi\endinsert}
\def\midspace#1{\gdef\thespace@{#1}\ifvmode\def\next
    {\futurelet\next\midspace@}\else\def\next{\nonvmodeerr@\midspace}\fi\next}
\def\midspace@{\ifx\next\Space@\def\next. {\futurelet\next\midspace@@}\else
     \def\next.{\futurelet\next\midspace@@}\fi\next.}
\def\midspace@@{\ifx\next\caption\let\next\midspace@@@\else
    \let\next\midspace@@@@\fi\next}
 \def\midspace@@@@{\midinsert\vbox to 
       \thespace@{}\endinsert}
\def\midspace@@@\caption#1{\midinsert\vbox to
    \thespace@{}\nobreak
      \smallskip
      \setbox\z@=\hbox{\noindent\ignorespaces#1\unskip}%
      \ifdim\wd\z@>\captionwidth@
    \centerline{\vbox{\hsize=\captionwidth@\noindent\ignorespaces#1\unskip}}%
    \else\centerline{\box\z@}\fi\endinsert}
\mathchardef\prime@="0230
\def\prime{{{}\prime@{}}}
\def\prim@s{\prime@\futurelet\next\pr@m@s}

\def\,{\relax\ifmmode\mskip\thinmuskip\else\thinspace\fi}
\def\!{\relax\ifmmode\mskip-\thinmuskip\else\negthinspace\fi}
\def\frac#1#2{{#1\over#2}}

\def\:{\nobreak\hskip.1111em{:}\hskip.3333em plus .0555em\relax}
\def\intic@{\mathchoice{\hskip5\p@}{\hskip4\p@}{\hskip4\p@}{\hskip4\p@}}
\def\negintic@
 {\mathchoice{\hskip-5\p@}{\hskip-4\p@}{\hskip-4\p@}{\hskip-4\p@}}
\def\intkern@{\mathchoice{\!\!\!}{\!\!}{\!\!}{\!\!}}
\def\intdots@{\mathchoice{\cdots}{{\cdotp}\mkern1.5mu
    {\cdotp}\mkern1.5mu{\cdotp}}{{\cdotp}\mkern1mu{\cdotp}\mkern1mu
      {\cdotp}}{{\cdotp}\mkern1mu{\cdotp}\mkern1mu{\cdotp}}}
\newcount\intno@             
\def\iint{\intno@=\tw@\futurelet\next\ints@} 
\def\iiint{\intno@=\thr@@\futurelet\next\ints@}
\def\iiiint{\intno@=4 \futurelet\next\ints@}
\def\idotsint{\intno@=\z@\futurelet\next\ints@}
\def\ints@{\findlimits@\ints@@}
\newif\iflimtoken@
\newif\iflimits@
\def\findlimits@{\limtoken@false\limits@false\ifx\next\limits
 \limtoken@true\limits@true\else\ifx\next\nolimits\limtoken@true\limits@false
    \fi\fi}
\def\multintlimits@{\intop\ifnum\intno@=\z@\intdots@
  \else\intkern@\fi
    \ifnum\intno@>\tw@\intop\intkern@\fi
     \ifnum\intno@>\thr@@\intop\intkern@\fi\intop}
\def\multint@{\int\ifnum\intno@=\z@\intdots@\else\intkern@\fi
   \ifnum\intno@>\tw@\int\intkern@\fi
    \ifnum\intno@>\thr@@\int\intkern@\fi\int}
\def\ints@@{\iflimtoken@\def\ints@@@{\iflimits@
   \negintic@\mathop{\intic@\multintlimits@}\limits\else
    \multint@\nolimits\fi\eat@}\else
     \def\ints@@@{\multint@\nolimits}\fi\ints@@@}
\def\Sb{_\bgroup\vspace@
        \baselineskip=\fontdimen10 \scriptfont\tw@
        \advance\baselineskip by \fontdimen12 \scriptfont\tw@
        \lineskip=\thr@@\fontdimen8 \scriptfont\thr@@
        \lineskiplimit=\thr@@\fontdimen8 \scriptfont\thr@@
        \Let@\vbox\bgroup\halign\bgroup \hfil$\scriptstyle
            {##}$\hfil\cr}
\def\endSb{\crcr\egroup\egroup\egroup}
\def\Sp{^\bgroup\vspace@
        \baselineskip=\fontdimen10 \scriptfont\tw@
        \advance\baselineskip by \fontdimen12 \scriptfont\tw@
        \lineskip=\thr@@\fontdimen8 \scriptfont\thr@@
        \lineskiplimit=\thr@@\fontdimen8 \scriptfont\thr@@
        \Let@\vbox\bgroup\halign\bgroup \hfil$\scriptstyle
            {##}$\hfil\cr}
\def\endSp{\crcr\egroup\egroup\egroup}
\def\Let@{\relax\iffalse{\fi\let\\=\cr\iffalse}\fi}
\def\vspace@{\def\vspace##1{\noalign{\vskip##1 }}}
\def\aligned{\,\vcenter\bgroup\vspace@\Let@\openup\jot\m@th\ialign
  \bgroup \strut\hfil$\displaystyle{##}$&$\displaystyle{{}##}$\hfil\crcr}
\def\endaligned{\crcr\egroup\egroup}
\def\matrix{\,\vcenter\bgroup\Let@\vspace@
    \normalbaselines
  \m@th\ialign\bgroup\hfil$##$\hfil&&\quad\hfil$##$\hfil\crcr
    \mathstrut\crcr\noalign{\kern-\baselineskip}}
\def\endmatrix{\crcr\mathstrut\crcr\noalign{\kern-\baselineskip}\egroup
                \egroup\,}
\newtoks\hashtoks@
\hashtoks@={#}
\def\format{\crcr\egroup\iffalse{\fi\ifnum`}=0 \fi\format@}
\def\format@#1\\{\def\preamble@{#1}%
  \def\c{\hfil$\the\hashtoks@$\hfil}%
  \def\r{\hfil$\the\hashtoks@$}%
  \def\l{$\the\hashtoks@$\hfil}%
  \setbox\z@=\hbox{\xdef\Preamble@{\preamble@}}\ifnum`{=0 \fi\iffalse}\fi
   \ialign\bgroup\span\Preamble@\crcr}

\def\cases{\left\{\,\vcenter\bgroup\vspace@
     \normalbaselines\openup\jot\m@th
       \Let@\ialign\bgroup$##$\hfil&\quad$##$\hfil\crcr
      \mathstrut\crcr\noalign{\kern-\baselineskip}}

\newif\iftagsleft@
\tagsleft@true
\def\TagsOnRight{\global\tagsleft@false}
\def\tag#1$${\iftagsleft@\leqno\else\eqno\fi
 \hbox{\def\pagebreak{\global\postdisplaypenalty-\@M}%
 \def\nopagebreak{\global\postdisplaypenalty\@M}\rm(#1\unskip)}%
  $$\postdisplaypenalty\z@\ignorespaces}
\interdisplaylinepenalty=\@M
\def\allowdisplaybreak@{\def\allowdisplaybreak{\noalign{\allowbreak}}}
\def\displaybreak@{\def\displaybreak{\noalign{\break}}}
\def\align#1\endalign{\def\tag{&}\vspace@\allowdisplaybreak@\displaybreak@
  \iftagsleft@\lalign@#1\endalign\else
   \ralign@#1\endalign\fi}
\def\ralign@#1\endalign{\displ@y\Let@\tabskip\centering\halign to\displaywidth
     {\hfil$\displaystyle{##}$\tabskip=\z@&$\displaystyle{{}##}$\hfil
       \tabskip=\centering&\llap{\hbox{(\rm##\unskip)}}\tabskip\z@\crcr
             #1\crcr}}
\def\lalign@
 #1\endalign{\displ@y\Let@\tabskip\centering\halign to \displaywidth
   {\hfil$\displaystyle{##}$\tabskip=\z@&$\displaystyle{{}##}$\hfil
   \tabskip=\centering&\kern-\displaywidth
        \rlap{\hbox{(\rm##\unskip)}}\tabskip=\displaywidth\crcr
               #1\crcr}}
\def\overrightarrow{\mathpalette\overrightarrow@}
\def\overrightarrow@#1#2{\vbox{\ialign{$##$\cr
    #1{-}\mkern-6mu\cleaders\hbox{$#1\mkern-2mu{-}\mkern-2mu$}\hfill
     \mkern-6mu{\to}\cr
     \noalign{\kern -1\p@\nointerlineskip}
     \hfil#1#2\hfil\cr}}}
\def\overleftarrow{\mathpalette\overleftarrow@}
\def\overleftarrow@#1#2{\vbox{\ialign{$##$\cr
     #1{\leftarrow}\mkern-6mu\cleaders\hbox{$#1\mkern-2mu{-}\mkern-2mu$}\hfill
      \mkern-6mu{-}\cr
     \noalign{\kern -1\p@\nointerlineskip}
     \hfil#1#2\hfil\cr}}}
\def\overleftrightarrow{\mathpalette\overleftrightarrow@}
\def\overleftrightarrow@#1#2{\vbox{\ialign{$##$\cr
     #1{\leftarrow}\mkern-6mu\cleaders\hbox{$#1\mkern-2mu{-}\mkern-2mu$}\hfill
       \mkern-6mu{\to}\cr
    \noalign{\kern -1\p@\nointerlineskip}
      \hfil#1#2\hfil\cr}}}
\def\underrightarrow{\mathpalette\underrightarrow@}
\def\underrightarrow@#1#2{\vtop{\ialign{$##$\cr
    \hfil#1#2\hfil\cr
     \noalign{\kern -1\p@\nointerlineskip}
    #1{-}\mkern-6mu\cleaders\hbox{$#1\mkern-2mu{-}\mkern-2mu$}\hfill
     \mkern-6mu{\to}\cr}}}
\def\underleftarrow{\mathpalette\underleftarrow@}
\def\underleftarrow@#1#2{\vtop{\ialign{$##$\cr
     \hfil#1#2\hfil\cr
     \noalign{\kern -1\p@\nointerlineskip}
     #1{\leftarrow}\mkern-6mu\cleaders\hbox{$#1\mkern-2mu{-}\mkern-2mu$}\hfill
      \mkern-6mu{-}\cr}}}
\def\underleftrightarrow{\mathpalette\underleftrightarrow@}
\def\underleftrightarrow@#1#2{\vtop{\ialign{$##$\cr
      \hfil#1#2\hfil\cr
    \noalign{\kern -1\p@\nointerlineskip}
     #1{\leftarrow}\mkern-6mu\cleaders\hbox{$#1\mkern-2mu{-}\mkern-2mu$}\hfill
       \mkern-6mu{\to}\cr}}}
\def\sqrt#1{\radical"270370 {#1}}
\def\dots{\relax\ifmmode\let\next=\ldots\else\let\next=\tdots@\fi\next}
\def\tdots@{\unskip\ \tdots@@}
\def\tdots@@{\futurelet\next\tdots@@@}
\def\tdots@@@{$\mathinner{\ldotp\ldotp\ldotp}\,
   \ifx\next,$\else
   \ifx\next.\,$\else
   \ifx\next;\,$\else
   \ifx\next:\,$\else
   \ifx\next?\,$\else
   \ifx\next!\,$\else
   $ \fi\fi\fi\fi\fi\fi}
\def\text{\relax\ifmmode\let\next=\text@\else\let\next=\text@@\fi\next}
\def\text@@#1{\hbox{#1}}
\def\text@#1{\mathchoice
 {\hbox{\everymath{\displaystyle}\def\textfonti{\the\textfont1 }%
    \def\textfontii{\the\textfont2 }\textdef@@ T#1}}
 {\hbox{\everymath{\textstyle}\def\textfonti{\the\textfont1 }%
    \def\textfontii{\the\textfont2 }\textdef@@ T#1}}
 {\hbox{\everymath{\scriptstyle}\def\textfonti{\the\scriptfont1 }%
   \def\textfontii{\the\scriptfont2 }\textdef@@ S\rm#1}}
 {\hbox{\everymath{\scriptscriptstyle}\def\textfonti{\the\scriptscriptfont1 }%
   \def\textfontii{\the\scriptscriptfont2 }\textdef@@ s\rm#1}}}
\def\textdef@@#1{\textdef@#1\rm \textdef@#1\bf
   \textdef@#1\sl \textdef@#1\it}

\def\textdef@#1#2{\def\next{\csname\expandafter\eat@\string#2fam\endcsname}%
\if S#1\edef#2{\the\scriptfont\next\relax}%
 \else\if s#1\edef#2{\the\scriptscriptfont\next\relax}%
 \else\edef#2{\the\textfont\next\relax}\fi\fi}
\scriptfont\itfam=\tenit \scriptscriptfont\itfam=\tenit
\scriptfont\slfam=\tensl \scriptscriptfont\slfam=\tensl
\mathcode`\0="0030
\mathcode`\1="0031
\mathcode`\2="0032
\mathcode`\3="0033
\mathcode`\4="0034
\mathcode`\5="0035
\mathcode`\6="0036
\mathcode`\7="0037
\mathcode`\8="0038
\mathcode`\9="0039
\def\Cal{\relax\ifmmode\let\next=\Cal@\else
     \def\next{\errmessage{Use \string\Cal\space only in math mode}}\fi\next}
\def\Cal@#1{{\fam2 #1}}
\def\bold{\relax\ifmmode\let\next=\bold@\else
   \def\next{\errmessage{Use \string\bold\space only in math
      mode}}\fi\next}\def\bold@#1{{\fam\bffam #1}}
\mathchardef\Gamma="0000
\mathchardef\Delta="0001
\mathchardef\Theta="0002
\mathchardef\Lambda="0003
\mathchardef\Xi="0004
\mathchardef\Pi="0005
\mathchardef\Sigma="0006
\mathchardef\Upsilon="0007
\mathchardef\Phi="0008
\mathchardef\Psi="0009
\mathchardef\Omega="000A
\mathchardef\varGamma="0100
\mathchardef\varDelta="0101
\mathchardef\varTheta="0102
\mathchardef\varLambda="0103
\mathchardef\varXi="0104
\mathchardef\varPi="0105
\mathchardef\varSigma="0106
\mathchardef\varUpsilon="0107
\mathchardef\varPhi="0108
\mathchardef\varPsi="0109
\mathchardef\varOmega="010A
\font\dummyft@=dummy
\fontdimen1 \dummyft@=\z@
\fontdimen2 \dummyft@=\z@
\fontdimen3 \dummyft@=\z@
\fontdimen4 \dummyft@=\z@
\fontdimen5 \dummyft@=\z@
\fontdimen6 \dummyft@=\z@
\fontdimen7 \dummyft@=\z@
\fontdimen8 \dummyft@=\z@
\fontdimen9 \dummyft@=\z@
\fontdimen10 \dummyft@=\z@
\fontdimen11 \dummyft@=\z@
\fontdimen12 \dummyft@=\z@
\fontdimen13 \dummyft@=\z@
\fontdimen14 \dummyft@=\z@
\fontdimen15 \dummyft@=\z@
\fontdimen16 \dummyft@=\z@
\fontdimen17 \dummyft@=\z@
\fontdimen18 \dummyft@=\z@
\fontdimen19 \dummyft@=\z@
\fontdimen20 \dummyft@=\z@
\fontdimen21 \dummyft@=\z@
\fontdimen22 \dummyft@=\z@
\def\fontlist@{\\{\tenrm}\\{\sevenrm}\\{\fiverm}\\{\teni}\\{\seveni}%
 \\{\fivei}\\{\tensy}\\{\sevensy}\\{\fivesy}\\{\tenex}\\{\tenbf}\\{\sevenbf}%
 \\{\fivebf}\\{\tensl}\\{\tenit}\\{\tensmc}}
\def\dodummy@{{\def\\##1{\global\let##1=\dummyft@}\fontlist@}}
\newif\ifsyntax@
\newcount\countxviii@
\def\newtoks@{\alloc@5\toks\toksdef\@cclvi}
\def\nopages@{\output={\setbox\z@=\box\@cclv \deadcycles=\z@}\newtoks@\output}
\def\syntax{\syntax@true\dodummy@\countxviii@=\count18
\loop \ifnum\countxviii@ > \z@ \textfont\countxviii@=\dummyft@
   \scriptfont\countxviii@=\dummyft@ \scriptscriptfont\countxviii@=\dummyft@
     \advance\countxviii@ by-\@ne\repeat
\dummyft@\tracinglostchars=\z@
  \nopages@\frenchspacing\hbadness=\@M}
\def\magstep#1{\ifcase#1 1000\or
 1200\or 1440\or 1728\or 2074\or 2488\or 
 \errmessage{\string\magstep\space only works up to 5}\fi\relax}
{\lccode`\2=`\p \lccode`\3=`\t 
 \lowercase{\gdef\tru@#123{#1truept}}}

\def\scaletype#1{\mag=#1\relax
 \hsize=\expandafter\tru@\the\hsize
 \vsize=\expandafter\tru@\the\vsize
 \dimen\footins=\expandafter\tru@\the\dimen\footins}

\def\scalefont#1#2\andcallit#3{\edef\font@{\the\font}#1\font#3=
  \fontname\font\space scaled #2\relax\font@}
\def\Mag@#1#2{\ifdim#1<1pt\multiply#1 #2\relax\divide#1 1000 \else
  \ifdim#1<10pt\divide#1 10 \multiply#1 #2\relax\divide#1 100\else
  \divide#1 100 \multiply#1 #2\relax\divide#1 10 \fi\fi}
\def\scalelinespacing#1{\Mag@\baselineskip{#1}\Mag@\lineskip{#1}%
  \Mag@\lineskiplimit{#1}}
\def\wlog#1{\immediate\write-1{#1}}
\catcode`\@=\active

\magnification=\magstep1

\centerline{\bf Ordinal Recursion
 Theory}
\vskip.3in
\centerline{C. T. Chong}

\centerline{\it National University of Singapore}
\vskip.2in
\centerline{S. D. Friedman$^1$}

\centerline{\it Massachusetts Institute
of Technology}
\footnote{Preparation of this
paper was supported by NSF Grant \#9205530}

\vskip.4in
\centerline{\bf Introduction}
\vskip.3in

In a fundamental paper, Kreisel and
Sacks [1965] initiated the study of
``metarecursion theory'', an analog of
classical recursion theory where
$\omega$ is replaced by Church-Kleene
$\omega_1$, the least non-recursive
ordinal. Subsequently, Sacks and his
school developed recursion theory on
arbitrary $\Sigma_1$-admissible
ordinals,  now known as
``$\alpha$-recursion theory''. 

In
Section I of the present article, we
present the basic concepts and
techniques of this theory, putting
particular emphasis on the main new
ideas that have been introduced to
study recursion-theoretic problems
assuming only
$\Sigma_1$-admissibility on a domain
greater than $\omega$.
As $\Sigma_1$-admissibility is easily
lost under relativization, we turn to 
``$\beta$-recursion theory''  (Section
II) which attempts to develop recursion
theory on arbitrary limit ordinals. 
In Section III,  the final part of this article, we
take up the topic of ``admissibility
spectra'', where instead of studying
the definability of subsets of a fixed
$\Sigma_1$-admissible ordinal, we ask:
given a set $X$, which are the ordinals
$\Sigma_1$-admissible relative to $X$ ?

The reader will notice that 
Jensen's work on the fine structure
theory of G\"odel's $L$ features
prominently throughout.
Indeed a major development of ordinal
recursion theory is the infusion of
set-theoretic ideas in studying
recursion-theoretic problems. The
unmistakeable presence of a  strong
set-theoretic flavor in   the subject of admissibility
spectra is especially pronounced. We
thus view the appearance of Jensen's
paper (Jensen [1972], preliminary
copies of which had been circulated
earlier), at a time when ordinal
recursion theory was being developed, to
be a fortuitous happening.  

Some of
the techniques and ideas which
were invented in ordinal recursion
theory have recently found applications in
``recursion
theory on fragments of Peano
arithmetic''. This is an unexpected
turn of events which signal a basic
unity among various fields in
recursion theory and fine structure
theory. 
We  touch briefly on this work
at the end of Section I.
\vskip.3in
\centerline{\bf I. $\alpha$-Recursion Theory}
\vskip.3in
\noindent

\vskip.3in
We begin with some basic
 notions. 
Recall  G\"odel's
constructible universe $L$, defined as $\cup
\{L_\alpha|\alpha \ \text{an ordinal}\}$.
   A limit   ordinal $\alpha$
is $\Sigma_n$-{\it admissible} 
 if
$L_\alpha$ satisfies the
replacement axiom for $\Sigma_n$
formulas (with parameters in $L_\alpha$) in 
ZF set theory.  If $\alpha$ is
$\Sigma_n$-admissible for some $n\ge 1$,
there is a $\Sigma_1(L_\alpha)$
bijection between $\alpha$ and
$L_\alpha$, allowing one to identify
these  two objects if and when necessary.
$\Sigma_1$-admissible ordinals are
sometimes referred to  simply as
admissible  ordinals.
 Unless otherwise specified, we fix
$\alpha$ to be an  admissible
ordinal henceforth.

A set $K\subset \alpha$ is
$\alpha$-{\it finite} if $K\in
L_\alpha$. A function is partial $\alpha$-{\it
recursive} if its graph is $\Sigma_1
(L_\alpha)$. A set   is $\alpha$-{\it
recursively enumerable} ($\alpha$-RE)
 if it is the
domain of a partial $\alpha$-recursive
function. $A\subset \alpha$ is $\alpha$-recursive if 
both $A$ and $\alpha\setminus A$ are
$\alpha$-RE. In terms of definability, a set is
$\alpha$-recursive if and only if it is
$\Delta_1(L_\alpha)$. It is
$\alpha$-finite if and only if it is
$\alpha$-recursive and bounded in $\alpha$.

All the basic results in classical
recursion theory, for example those
covered in the first seven chapters of
Rogers [1967], hold for all $\Sigma_1$
admissible ordinals.  Thus
a set $K\subset  \alpha$ is RE~if and
only if it is the range of  a total
$\alpha$-recursive function; there is
an effective  
(i.e.~$\Sigma_1(L_\alpha)$ definable)
enumeration of all $\alpha$-finite
sets and all $\alpha$-RE sets;
Kleene's Recursion Theorem is true for
each $\alpha$. We denote by $W_e$ the
$e$th $\alpha$-RE~set and by $K_e$
the $e$th $\alpha$-finite set under the
respective effective enumerations.

\vskip.2in
\noindent{\bf Reducibility} The notion of
reducibility provides a means of
comparing the relative complexity of
subsets of $\alpha$. Given  $A\subset
\alpha$, define by the collection of
neighborhood conditions of $A$ the set
$$
N(A)=\{(c,d)|K_c\subset A\ \& \
K_d\subset \alpha\setminus A\}.
$$
We say that $A$ is $\alpha$-RE in
$B\subset \alpha$ if there is an $e$
such that for all $x<\alpha$,
$$
x\in A\leftrightarrow (\exists
c)(\exists d)[(x,c,d)\in W_e\ \& \
(c,d)\in N(B)\}.
$$
$A$ is {\it weakly} $\alpha$-{\it
recursive in} $B$, written
$A\le_{w\alpha} B$,  if $A$ and $\bar{A}$
are $\alpha$-RE in $B$. 
Define by $A^*$ the set $\{u|K_u\subset
A\}$. Then $A$ is $\alpha$-{\it
recursive in} $B$, written $A\le_\alpha
B$,  if $A^*$ and
$\bar{A^*}$ are $\alpha$-RE in $B$.

Thus $A$ is $\alpha$-recursive in $B$
provided there is an algorithm such that
for any given $\alpha$-finite set $K$, it
is possible to use the algorithm, with
$B$ as an oracle, to conclude within
$\alpha$-finite time whether $K$ is a
subset of $A$ or 
disjoint from $A$.
It is not difficult to verify that
$\le_\alpha$ is reflexive and
transitive. $A$ and $B$ are said to
have the same $\alpha$-degree, written
$A\equiv_\alpha B$, if
$A\le_\alpha B$ and $B\le_\alpha A$. 
 Although
$\le_{w\omega}$ is equivalent to
$\le_\omega$ and therefore transitive, 
$\le_{w\alpha}$ is not transitive in general.

The least complicated $\alpha$-degree,
which we denote by {\bf 0},
is the $\alpha$-recursive degree, which
consists of $\alpha$-recursive sets. A
degree is an $\alpha$-RE~degree if it
contains an $\alpha$-RE~set.
There is a {\it greatest} $\alpha$-RE~degree
{\bf 0$^\prime$} which contains the
$\alpha$-RE~set
$\emptyset^\prime=\{(x,e)|x\in W_e\}$,
in which every $\alpha$-RE~set is $\alpha$-recursive.

There is an analog of  Church's  thesis for
$\alpha$-recursion theory  which we
shall appeal to in this article.
 This thesis allows a more
informal presentation  of the topics to
be covered, emphasizing intuition over formalism. 

 The key motivation of ordinal
recursion   lies in  the search
for a ``generalized'' recursion theory.
It is evident that the  notion of
effective   computation applies to a wider class of
mathematical structures, as exemplified
 in  Kleene's work on ordinal
notations (Kleene [1938]).
Kreisel and Sacks [1965] initiated the study
of recursion theory on Church-Kleene
$\omega_1$, and this led to the
subsequent introduction of the theory of
$\Sigma_1$ admissible ordinals by Sacks
and his school.  

A  closer examination reveals  that if
the
full replacement axiom is assumed in
$L_\alpha$, then 
 many
difficult proofs in classical recursion
theory go through almost routinely,
without major  modifications of the
classical construction
(there are exceptions: cf.~the section
on maximal sets).
>From the point of view of effective
computation, where ``$\Sigma_1$''-ness is
 identified with ``effectively enumerable'', it
should be sufficient 
 to assume only
$\Sigma_1$-replacement axiom to arrive
at a satisfactory recursion theory (though  all is
not lost even when
this crucial assumption is removed in
$\beta$-recursion theory, see Section II
).  This view is  supported by the
successful solution of  Post's problem for all
admissible ordinals 
 (Theorem
2 below).

We give here some examples of admissible ordinals:

\itemitem{(a)} $\alpha=\omega$, the
classical case.  Then $\alpha$ is
$\Sigma_n$-admissible for all
$n<\omega$. The same conclusion holds
for any  regular constructible cardinal;

\itemitem{(b)} $\alpha=\omega_1^{CK}$,
Church-Kleene $\omega_1$.
Here $\alpha$ is $\Sigma_1$ but not
$\Sigma_2$-admissible. There is also a
$\Sigma_1(L_\alpha)$ map from $\alpha$
into $\omega$ ($\omega$ is called the
$\Sigma_1$-projectum of $\alpha$).
A  subset of $\omega$ is
$\alpha$-RE~if and only if it is
$\Pi_1^1$ definable;

\itemitem{(c)} $\alpha=\delta_2^1$,  the least
ordinal which is not the order type of
a $\Delta_2^1$ set of natural numbers.
In this case a
set of natural numbers is
$\alpha$-RE~if and only if it is
$\Sigma^1_2$ definable.

\itemitem{(d)}
$\alpha=\aleph^L_\omega$, the
$\omega$th constructible cardinal of
$L$. There is a $\Sigma_2(L_\alpha)$
cofinal  map from $\omega$ into
$\alpha$. Hence $\alpha$ is not
$\Sigma_2$-admissble. On the other
hand, every infinite cardinal in a well-founded
model of ZF is $\Sigma_1$-admissible.
\vskip.3in
\noindent{\bf  Regularity} 
\vskip.2in
A set $A\subset \alpha$ is {\it
regular} (in $\alpha$) if its restriction to every
$\gamma <\alpha$ is $\alpha$-finite.
It follows  that every set of natural numbers  is
regular. On the other hand, in
$\omega_1^{CK}$, Kleene's $\Cal O$, a
complete $\Pi_1^1$-set of natural
numbers, is not regular (in
$\omega_1^{CK}$), even though it is
bounded and  $\omega_1^{CK}$-RE.
Non-regularity  is a major feature which
distinguishes 
ordinal recursion theory from the
classical theory. Non-regular
$\alpha$-RE sets are sets with bounded
parts which cannot be enumerated in
$\alpha$-finite time. Their existence
renders some of the standard techniques
ineffective. Nevertheless, at
least for $\alpha$-RE sets and for
the study of $\alpha$-degrees, this
difficulty can be circumvented:

\vskip.3in
\noindent {\bf   Theorem 1 (Sacks [1966])}
{\sl Let $\alpha$ be
admissible. Then every
$\alpha$-RE~degree contains a 
regular $\alpha$-RE~set.}
   \vskip.2in
Maass [1978]
 showed that there is a
parameter free $\Sigma_1$ function $f$ such
that for any $\alpha$ and $e<\alpha$,
$W_e$ and $W_{f(e)}$ have the same
$\alpha$-degree and $W_{f(e)}$ is regular.

\vskip.3in
\noindent
{\bf Definability}
\vskip.2in
Jensen's work on the fine structure of
$L$ [1972] turns out to be a key
 component in  the development of 
ordinal recursion
theory, a development which  arguably
exemplifies the successful integration of
set-theoretic and recursion-theoretic
ideas. In retrospect, the secret to the
solutions of such basic problems  as
Post's problem, Sacks Splitting
Theorem, and the Density Theorem for
all admissible $\alpha$,
rests on the insight that 
the complexities of the classical
constructions, with the
intervention of 
 fine structure theory,   may be refined to achieve
the goals,  provided that
one chooses the appropriate definable
objects within $\alpha$ to carry out
the necessary priority arguments.
On the other hand, in certain cases where
such an approach fails,
it is shown that 
 the problems being considered have
negative solutions.
 Problems such as the
existence of maximal sets, and ordering
of
$\alpha$-degrees above {\bf 0$^\prime$},
are examples.

We list here several important objects
in fine structure theory that play
pivotal roles in ordinal recursion. Let
$B\subset \alpha$.

The $\Sigma_n$-cofinality of
$(L_\alpha,B)$ is defined to be the least
$\gamma$ for which there is a
$\Sigma_n(L_\alpha,B)$ function from
$\gamma$ cofinally into $\alpha$.
We denote this ordinal by $\kappa_n(B)$,
or simply write it as
$\Sigma_n$-cofinality $(\alpha, B)$.
Clearly $\kappa_n(\emptyset)=\alpha$ if and only
if $\alpha$ is $\Sigma_n$-admissible.

The $\Sigma_n^B$-projectum of $\alpha$,
denoted $\alpha^*_n(B)$ (or sometimes
$\Sigma_n$-projectum $(\alpha, B)$), is
the least ordinal $\gamma\le\alpha$ for which there is a
$\Sigma_n(L_\alpha,B)$ map from $\alpha$ into
 $\gamma$. If $B=\emptyset$, Jensen's
 theory  provides
several characterizations of this
ordinal: (a)
it is the largest limit ordinal less
than or equal to $\alpha$ in which
every bounded $\Sigma_n(L_\alpha)$ set
is $\alpha$-finite; (b) it is the least
ordinal $\gamma$ for which there is a
$\Sigma_n(L_\alpha)$ map from a subset
of $\gamma$ onto  $\alpha$.
 When $B$ is
$\alpha$-RE and regular,  the
above characterization continues to hold with
$\Sigma_n(L_\alpha)$ replaced by
$\Sigma_n(L_\alpha, B)$, given and used
in Shore's proof of the Density Theorem [1976].

We use the notations  $\alpha^*_n$ and $\kappa_n$
when $B$ is empty. When $n=1$, we omit
the subscript $1$ and
simply write $\alpha^*$ and $\alpha^*(B)$
instead.  In $\alpha$-recursion
theory, it is important to  present
the set of requirements with as short
a 
list as possible. The ordinal
$\alpha^*$ or $\alpha^*(B)$ are often
used for this purpose.

We say that $\lambda<\alpha$  is an
$\alpha$-cardinal if
 there is no
$\alpha$-finite injection of $\lambda$
into a smaller  ordinal. 
 If $\alpha^*_n <\alpha$ (or if $
\kappa_n(B)<\alpha$),
it is not difficult to prove that it is
an $\alpha$-cardinal. And 
$\alpha^*<\alpha$ implies that it
is the greatest $\alpha$-cardinal.

A set  $B$ is {\it
hyperregular} if $\kappa_1(B)=\alpha$.
In other words, $B$ is
hyperregular if $(L_\alpha , B)$ is a
$\Sigma_1$-admissible structure. It is
not difficult to verify that every
$\alpha$-recursive set is 
hyperregular.
However, there are $\alpha$-RE~sets
which do not satisfy hyperregularity.
As an example, consider the set $B$ of
non-cardinals in $\alpha=\aleph_\omega^L$.
This is an $\alpha$-RE~set 
whose
complement is of order type $\omega$.
Then $f(n)=n$th member of
$\alpha\setminus B$ is a function weakly
$\alpha$-recursive in $B$, mapping
$\omega$ unboundedly into $\alpha$.
It turns out that for
$\alpha=\aleph^L_\omega$, the only
nonhyperregular $\alpha$-RE~set is of
complete degree {\bf 0$^\prime$}, and
every set of $\alpha$-degree above this
is non-hyperregular, while any set
which does not compute $B$ defined above
is hyperregular (under the axiom of
constructibility). In particular, every
incomplete $\alpha$-RE~set is hyperregular.

Hyperregularity is a strong condition
which ensures that computations carried
out on $\alpha$-finite sets using
oracles are completed in
$\alpha$-finite steps. Its
recursion-theoretic consequences are
significant: For example, for the 
$\alpha$ considered above, there is no
incomplete $\alpha$-RE~set whose ``jump'' (an analog
of the classical notion) is strictly
above $\emptyset^\prime$ (hence no
incomplete ``high set'') (Shore [1976a]). On
the other hand, many tools, such as
priority arguments, are not
relatiivizable  
to non-hyperregular sets. This
introduces additional complications to
the study of $\alpha$-recursion theory.
Different techniques are needed in many
cases and, in the most extreme case,
non-hyperregularity leads to radically different
degree-theoretic results (see Section 5
below).

\vskip.3in 
\noindent{\bf The $\alpha$-Finite Injury
Priority Method}
\vskip.2in
The method of finite injury priority argument
introduced by Friedberg  and
Muchnik  to solve Post's  problem marked the advent of
modern recursion theory. 
This
 method has since been joined by a
variety of highly complex and ingenious 
techniques invented to handle
problems about RE~sets and their degrees,
of which the
Friedberg-Muchnik proof is now
seen to be the simplest.  Solving
Post's problem may indeed be considered
the  first important test for   any
reasonable ordinal
recursion theory. 

\vskip.3in
\noindent
{\bf  
 Theorem 2  (Friedberg-Muchnik
Theorem)}  {\sl Let $\alpha$ be an
admissible ordinal. There
exist $\alpha$-RE~sets $A$ and $B$
with
incomparable $\alpha$-degrees.}

\vskip.3in
\noindent{\bf  Corollary (Solution of
Post's Problem)}  {\sl  There exists an
incomplete, $\alpha$-RE and non-$\alpha$-recursive degree}.
\vskip.3in

We sketch the solution  by Sacks
and Simpson [1972].  This is a   proof
which has a strong model-theoretic and
set-theoretic  flavor, in contrast to
that of  Lerman
 [1972] which  has a
stronger recursion-theoretic tilt  (it is
worth noting that Lerman's approach may be
refined to provide a parameter   free
construction of the sets $A$ and $B$,
yielding a uniform  solution, for all
admissible ordinals, to
Post's problem (Lerman, unpublished)).

\vskip.3in

\noindent {\bf Proof of Theorem
1:} 
Consider requirements of
the type
$$
R_e: \{e\}^A\ne B
$$
and those with the roles of $A$ and
$B$ reversed. The basic strategy is
 to  diagonalize against equalities whenever possible
while preserving   computations,
respectingg   requirements of higher
priority if and when necessary.
 The strategy succeeds in the classical
theory because
   (a)  for any $e_0$,
there is a stage after which  no requirement
of higher priority than $e_0$ gets
injured, and (b) it can be established that each
requirement $R_e$  gets imjured at most
finitely many  (indeed $2^e$) times.   Closer inspection shows
that (a) is  essentially a $\Sigma_2$
condition, and when satisfied,  is
sufficient to  derive (b).
What the construction demands then is for
a $\Sigma_1$-admissible ordinal $\alpha$
 to
perform a $\Sigma_2$ task.
(We will  see later that with the Density Theorem,
the required task is even more
onerous---at almost the
$\Sigma_3$ level).
\vskip.2in
Since $\alpha$ is in general not
$\Sigma_2$-admissible, a straighforward
adaptation of the original approach
will clearly fail. Instead, the
following 
two lemmas provide an insight into how
the difficulties may be overcome:

\vskip.2in
\noindent
{\bf    Sacks-Simpson Lemma} {\sl  Let
$\kappa$ be  a regular
$\alpha$-cardinal.
Suppose that $K$ is an $\alpha$-finite
set of $\alpha$-cardnality less than $\kappa$ such that $\{I_d|d\in K\}$  is a
simultaneous  $\alpha$-RE sequence of
 $\alpha$-finite sets each of which
has $\alpha$-cardinallity less
than  $\kappa$. Then $\cup_{d\in K}
I_{d}$  is  $\alpha$-finite and of
$\alpha$-cardinality less than    $\kappa$.}

\noindent

\vskip.2in
An ordinal $\sigma<\alpha$ is said to
be $\alpha$-stable if  $L_\sigma$ is a
$\Sigma_1$-elementary substructure of $L_\alpha$.

\vskip.2in
\noindent
{\bf  $\alpha$-Stability Lemma}
{\sl  If  $\omega<\alpha=\alpha^*$,  then it
is   a limit  of  $\alpha$-stable ordinals.}

\vskip.2in
The first step to the solution is to
provide a short indexing of
requirements with its associated list
of priorities. 
There are several cases to consider. First
suppose that $\alpha$ has a
greatest $\alpha$-cardinal $\kappa$.

 (a) $\kappa=\alpha^*<\alpha$.  In this case we use
the $\Sigma_1$-projectum $\alpha^*$ of
$\alpha$  to provide a  list of  requirements. Let $p$
be a one-one  $\alpha$-recursive
map from $\alpha$ into  $\alpha^*$. Requirement
$R_d$ is said to have higher priority
than requirement $R_e$ if $p(d)<p(e)$.
This  shorter list of indices ensures
that every $\alpha$-RE~set bounded in
$\alpha^*$ is $\alpha$-finite, an
essential feature that is needed during
the inductive stage to verify that
every requirement is satisfied.

Now commence with the construction
using the revised indexing of requirements.
For each  $e<\alpha$, let $I_e$
denote the  injury  set (defined in the
usual sense) associated
with $R_e$. The main observation here
is that if there is a regular
$\alpha$-cardinal $\rho\le \alpha^*$
such that $p(e)<\rho$ and $p(d)<p(e)$
implies that  $I_d$ has $\alpha$-cardinality less
than $\rho$,
 then the Sacks-Simpson Lemma
ensures that 
$I_e$ has $\alpha$-cardinality less than $\rho$
as well.  This is sufficient to show
that the requirement with the highest
priority after that $R_e$ is injured
less than $\rho$ times. Induction
hypothesis then allows one to conclude that
 every requirement is eventually satisfied.

\vskip.2in
By the $\alpha$-Stabiulity Lemma 2.4, let $\beta$ be the order
type of $\alpha$-stable ordinals above $\kappa$.

\vskip.2in
(b) $\alpha^*=\alpha$.  If $\beta=\alpha$
then we use the identity function for
priority listing, and modify the
classical construction slightly.
 Lemma 2.4  and the construction provides the necessary
tool to argue that every requirement
settles down before the next
$\alpha$-stable ordinal. 
If $\beta<\alpha$, then
 there is a $\Sigma_2(L_\alpha)$ bijection $p$
between  $\alpha$ and $\kappa\cdot\beta$
(since the property of ``being
$\alpha$-stable'' is $\Pi_1(L_\alpha)$).
We use $\kappa\cdot\beta$ to index the
requirements   and say that  $R_d$ has
higher priority than $R_e$ if $p(d)<p(e)$.
The  positions  of the  priorities
are given by  an $\alpha$-recursive
approximation $p^\prime$ of $p$.
This gives meaning to  ``the priority of  $R_e$
at stage $\sigma$ is $\nu$''. The
``final priority'' of $R_e$ is then
$p(e)$, which is the limit of $p^\prime
(\sigma,e)$ as $\sigma$ tends to $\alpha$.

Construction proceeds as before, using
$p^\prime$ to guide the priority
ordering at each stage.
The rules  governing the
injury of requirements in order of
priority at each stage are observed. 
Exploiting the
property of $\Sigma_1$-stability,
coupled with Lemma 2.4, ensures that
all requirements of priority at least
$\kappa\cdot \nu$ settle down by the
$\nu +1$-th $\alpha$-stable ordinal.

Finally, 
if  $\alpha$ is a limit
of $\alpha$-cardinals (analogous  to
$\omega$), one uses an
indexing provided by the identity
function  on $\alpha$. The argument
then proceeds as in Case (a).
\vskip.3in
\noindent
{\bf  The Density Theorem}
\vskip.2in
\noindent{\bf  Theorem 3
(Shore [1976]).} {\sl Let
{\bf b $<$ c} be $\alpha$-RE degrees.
Then there is an $\alpha$-RE degree
{\bf a} such that {\bf b $<$ a
$< $ c.}}
\vskip.2in
This theorem is one of the first
successful liftings of infinite injury
priority argument to ordinal recursion
theory. We sketch the key ideas here. 
Fix $B<_\alpha C$
to be regular $\alpha$-RE sets (Theorem
1). An $\alpha$-RE set $A$ of
intermediate $\alpha$-degree is to be constructed. 
\vskip.2in
\noindent{\bf  (Shore Incompleteness
Lemma)} {\sl 
Supposen $B$ is an incomplete $\alpha$-RE~set.
Then $\kappa_1(B)
\ge\alpha^*(B)$. 
Furthermore there is a
$\Sigma_1(L_\alpha, B)$ map from
$\kappa^*_1(B)$ onto $\alpha$.}

\vskip.2in
\noindent {\bf Proof.} We sketch the
proof of the first half of the lemma. Assume 
$\kappa^*_1(B) <\alpha^*(B)$. Let $D$ be a
regular $\alpha$-RE~set. We show that
$D\le_\alpha B$. Fix
$g:\kappa_1(B)\rightarrow \alpha$ to
be cofinal. Consider
$$
K=\{(\gamma,\delta)\vert
D\cap g(\gamma)\subset D^{g(\delta)}\}.
$$
Now $K$ is a $\Pi_1(B)$ set bounded
below $\alpha^*(B)$, and so is
$\alpha$-finite. Using it as a
parameter set, we see that $D\le_\alpha
B$. Hence $B$ is complete.
\vskip.2in
The Lemma  says essentially that if $B$
is incomplete, then $(L_\alpha, B)$ is
a {\it weakly admissible} structure. 
Weak admissibility allows many
$\Sigma_2(B)$  constructions, with
suitable modifications, to go
through (for example  Post's problem in
$\beta$-recursion theory,
cf.~Section II).

There are essentially three key
ingredients used in the proof of
Theorem 3: The
use of $\alpha^*(B)$
 for
a sufficiently short  listing of
the set of requirements; the
exploitation of the  blocking technique,
in which
requirements are grouped into
$\kappa_2(B)$ many
blocks of the same priority;
 and the use
of  $\kappa_1(B)$
and its associated cofinal function to
measure lengths of agreements between
computations in the course of the
construction. We elaborate the points below.

\vskip.2in

\itemitem{(a)}
 Let $p\le_{w\alpha} B$ be an injection
from $\alpha$ into $\alpha^*(B)$.
There is  a simultaneous
$\alpha$-recursive approximation
$\{p_\sigma\}$ of  $p$ such that for
all  $x$, $p_\sigma  (x)=p(x)$  for all
 sufficiently large  $\sigma$. 
 Requirements
are given a short list of length
$\alpha^*(B)$ using $p$. The principal feature of
this ordinal exploited in the proof of
the Density Theorem is that
every set $\alpha$-RE in $B$ and
bounded below $\alpha^*(B)$ is $\alpha$-finite.

\itemitem{(b)}  In the   construction there are altogether
$\kappa_2(B)$-many blocks of
requirements. Requirements in the same
block are accorded the same priority. This
reduces at once the number of injury sets
to
a manageable level. During 
verification step, one does induction
on $z<\kappa_2(B)$, and argues first of
all that the set of permanent injuries
inflicted on the computations is
bounded within each block, and secondly
that such a bound may be found in a
$\Sigma_2(L_\alpha, B)$ manner as a
function of $z$. The fact that
$z<\kappa_2(B)$ then ensures that a
uniform bound exists for all blocks
$z^\prime\le z$.

\itemitem{(c)}  $\kappa_1(B)$ is also
known as the {\it recursive cofinality
of\/} $B$.
Let $k: \kappa_1(B)\rightarrow \alpha$
be  a cofinal map weakly
$\alpha$-recursive in $B$.  There is a simultaneous
$\alpha$-recursive sequence of $\alpha$-recursive functions
$\{k_{\sigma}\}$  such  that  for each
$y<\kappa_1(B)$,
$k_\sigma|y=k|y$ for all
sufficiently  large  $\sigma$.
By the Shore Incompleteness Lemma  
 we may choose
 $k$
 to be a surjective map. The
calculations of lengths of agreement
between two computations will be based
on $k|y$, for $y<\kappa_1(B)$.
Furthermore, the surjectivity of $k$
allows the construction to pick up
every $\alpha$-finite set contained in 
$C$. During the
construction, such sets are coded
into $A$ (which in turn  causes complications).
 Each of these strategies
is designed to ensure that should $B$
be able to compute $A$, or $A$ compute
$C$, then in fact $C\le_\alpha B$, a contradiction.

\vskip.2in

There  are two  types of  requirements.
 The  positive requirements   $\{e\}^B\ne A$ for each  $e$, which  
attempts   to ensure  that  the  set $A$
to  be constructed  is not
$\alpha$-recursive  in  $B$,
   and
negative  requirements $\{e\}^A\ne C$
for each $e$, which  arranges
that $C$ is  not $\alpha$-recursive
in $A$.  These  requirements   are
grouped into  blocks indexed  by
$\kappa_2(B)$ with the aid of the
following lemma. Denote
$B^{<\sigma}$ to be the set of ordinals
enumerated in $B$ before stage
$\sigma$. Assume $\kappa_1(B)
>\omega$. Then
$B^{<\sigma}=B\cap \sigma$
(i.e.~$\sigma$ is $B$-correct) for
unboundedly many $\sigma$.
\vskip.2in
\noindent {\bf  Blocking Lemma}
{\sl There is a function
$g:\kappa_2(B)\rightarrow \alpha^*(B)$ which is
$\Sigma_2 (L_\alpha, B)$, together with
 a simultaneous
$\alpha$-recursive sequence
$\{g_\sigma\}$ of $g$ such that
\itemitem{(a)} $g_\sigma (z)\ge g(z)$ for all
sufficiently large $\sigma$;
\itemitem{(b)}
For all $z<\kappa_2(B)$,
and  for all
sufficiently large $B$-correct
$\sigma$, $g_\sigma|z= g|z$.}

\vskip.2in
We  say  that $\{g_\sigma\}$ is a {\it
tame   approximation} of  $g$ in view
of the Blocking Lemma. We shall only consider
$\kappa_1(B)>\omega$ here. The case when
$\kappa_1(B)=\omega$ is considerably
simpler. We say that a requirement with
index $e$
  is in block  $z$  if $p(e)
< g(z)$.
With the blocking lemma, it
makes sense using $\{g_\sigma\}$ to say that  a requirement is
``in block $z$ at stage
$\sigma$''. 
 Indeed by the Blocking Lemma, 
 if a
requirement is in block $z$, then it is
in block $z$ for all sufficiently large
$B$-correct $\sigma$. Furthermore, by tameness
property, this occurs uniformly in $z$,

As   in the  classical construction,
we code  the set $B$ into the  even
part of  $A$,  and  think  of the odd
part of $A$ as consisting of  triples
$(z,x,\sigma)$.  The three ordinals are
related via a length of agreement
function: Suppose
   at stage $\sigma$ 
there is an $e$ in  block $z$,   with
$A^{<\sigma}|k_\sigma(x)$ agreeing with 
$\{e\}^{B^{<\sigma}}_\sigma |k_\sigma(x)$.
Then ordinal $k_\sigma(x)$ is said to be the length of
agreement of computation for $e$ at stage $\sigma$.

This length   of  agreement  is
destroyed (i.e.~computation restarts)
at a later stage if new 
 elements below $\sigma$  enter either
 $A$ or $B$, since such occurances are
likely to invalidate any computations
reached so far. 
Those agreements which are never destroyed are
called permanent. They turn out to be
$\alpha$-recursively identifiable by
$B$.  Precaution is taken
so that if $(z,x,\sigma)$ enters $A$,
then  $K_{k_\sigma(x)}$, the
$k_\sigma(x)$th $\alpha$-finite set,
is contained in $C$. Since $C$ is
regular, this can be verified at some
stage. And since $k$ is a surjective
map, all the relevant $\alpha$-finite
sets will be  considered at some stage. The objective here is to code
enough of $C$ into $A$ so as to obtain
the following lemma:
\vskip.2in
\noindent {\bf  Lemma}
\itemitem{(i)} {\sl For each $e$ in block $z$,
$\{e\}^B\ne A$;}
\itemitem{(ii)} {\sl Within each block,
the permanent lengths of
agreemnt are bounded below $\alpha$.}

\noindent {\bf Proof:} The idea is that if (i)
fails, then using $B$ to identify   permanent
lengths of agreement,
one is able to compute $C$ (which are
coded in $A$) from $B$, 
 a contradiction.

A special feature of the blocking
technique is that requirements within
the same block work together to
achieve collectively  a
longer length of agreement. To prove
(ii), one uses the fact that the set $K$
of $e$'s in block $z$ for which 
$\{e\}^B$ is total on $k(x)$ and not
equal to $A|k(x)$  after an
agreement had been reached earlier, is a set
$\Sigma_1$ in $B$ and bounded below
$\alpha^*(B)$, hence $\alpha$-finite.
It is then sufficient to consider only
 $e\in z\setminus K$.
Repeating an
argument similar to that for (i) above on the
set $z\setminus K$, but
this time collectively on all the
computations that provide permanent lengths of agreement,
 shows that if  (ii) is false, then again
$C\le_\alpha B$.
\vskip.2in

Consider requirements $e$ in block $z$. A negative requirement is intended to
preserve computations of the form
$\{e\}^{A^{<\sigma}}|k_\sigma(x)$
to make it different from $C|k_\sigma(x)$. 
At stage $\sigma$, each requirement $e$ is assigned a marker
which is placed at the least ordinal
$\nu_{e,\sigma}$ greater than the negative facts used
about $A^{<\sigma}$ in the computation
above. The idea is that for as long as
markers stay, then no new ordinal below
their positions is allowed to enter
$A$. On the other hand, should a new element below
$\sigma$ enter $B$ at a later stage,
 then all markers assigned at stage
$\sigma$ are removed, clearing the way
for  ordinals below $\nu_{e,\sigma}$ to be added
to $A$ if and when necessary.
 These markers may reappear subsequently
(say at $\zeta>\sigma$) occupying different
positions  provided that, 
roughly speaking, there is  a  {\it
collective} length of
agreement between $C^{<\zeta}$ and
$\{\{e\}_\zeta^{A^{<\zeta}}\}$, $e$ in
block $z$, longer
than those achieved before.

The construction of the set $A$
involves the coding of $B$ into the
even part of $A$ (to ensure
$B\le_\alpha A$), and the manipulation
of positive and negative requirements.
A negative requirement (marker) is permanent if
it is never removed. 
The set of permanent negative
requirements within a block has to be
bounded else one argues that $C$ is
$\alpha$-recursive in $B$. Furthermore,
it can be arranged
 that within a block, the limit
inferior of the positions of
the  markers that stay behind at the
end of each stage of construction is
bounded below $\alpha$, and may be
computed from $C$. This allows
unboundedly many opportunities for
ordinals above certain level to enter
$A$, and is crucial to the success of
the construction. With this it is also 
possible to show that
 $C$ is not
$\alpha$-recursive in $A$. 

The final thread is to establish 
$A\le_\alpha C$. This is achieved
through arranging the construction so
that the set of permanent
negative requirements is
$\alpha$-recursive in $C$. We omit the details.

\vskip.3in
\noindent{\bf  Non-existence of
Maximal Sets}
\vskip.2in

In this and the next section, we give
two examples of 
problems which have negative solutions
in ordinal recursion. The first, due to
Lerman [1974], states roughly that
there is a lattice-theoretic property
of 
$\alpha$-RE~sets which is  inherently definably
countable. More precisely,
\vskip.2in
\noindent{\bf  Theorem 4} {\sl
There is a maximal $\alpha$-RE~set if
and only if there is a function $f$
which is $S_3$-definable mapping
$\alpha$ onto $\omega$.}

The notion of maximality is derived
from the classical one: 
 $M$ is maximal if and
only if its complement $\bar{M}$ is unbounded,
 and there is no
$\alpha$-RE~set which splits
$\bar{M}$ into two non-$\alpha$-finite
parts. We say that $f$ is $S_3$-definable if
there is an $\alpha$-recursive function
$f^\prime$ such that for all $x<\alpha$,
$$
\text{lim}_\tau\text{lim}_\sigma
f^\prime (\tau,\sigma,x)=f(x).
$$
Thus in our terminology, we may say
that there is a maximal
$\alpha$-RE set if and only if the
$S_3$-projectum of $\alpha$ is $\omega$.
This complete characterization of the
existence of maximal sets
 raises a very
interesting but apparetly quite difficult question: is there a
classification of recursion-theoretic
problems which  are inherently
linked to the cardinality of the
universe ? 

The following weak form of
Theorem 4
 shows how the size of
$\alpha$ has a bearing on the existence
of maximal sets:
\vskip.2in
\noindent{\bf  Theorem 5} {\sl If there is
a maximal $\alpha$-RE~set, then
$\alpha$ is countable.}
\vskip.2in
To prove this theorem, we consider
$\kappa_2$ which is
$\Sigma_2$-cofinality $(\alpha)$, and
$\alpha^*_2$, the $\Sigma_2$-projectum
of $\alpha$.

\vskip.2in
 
\noindent{\bf Proof of Theorem 5:} Let $M$ be
a maximal set.
We first claim
 that $\kappa_2\ge
\alpha^*_2$. To do this, build a
simultaneous $\alpha$-recursive
sequence of pairwise disjoint
$\alpha$-finite sets
$\{H_\nu\}_{\nu<\kappa_2}$ such that
$\bar{M}\cap (\cup_{\nu<\kappa_2}
H_\nu)$ is not $\alpha$-finite, and each $H_\nu$ contains at
most one member of $\bar{M}$. Now the set
$$
K=\{\nu|H_\nu\cap\bar{M}\ne\emptyset\}
$$
is a $\Sigma_2$ definable subset of
$\kappa_2$. If $\kappa_2<\alpha^*_2$,
then $K$ is $\alpha$-finite, in which
case it is possible to split $K$
into two non-empty parts $K_1$ and
$K_2$ so that $\cup_{\nu\in K_1} H_\nu$
and $\cup_{\nu\in K_2} H_\nu$ each
contains a non-$\alpha$-finite
unbounded subset
of $\bar{M}$, contradicting maximality
of $M$. Thus $\kappa_2\ge\alpha^*_2$.

Next we argue that $\kappa_2$ is in
fact countable. To do this, let
$\beta$ be the order type of $\bar{M}$.
 Partition $\alpha$
into an $\alpha$-recursive sequence of
pairwise disjoint $\alpha$-RE~sets $\{A_\nu\}_{\nu<\kappa_2}$.
Define
$$
B_\nu=\{\gamma|\exists
\sigma[\text{order type of}\ \gamma\setminus
M^\sigma]\in A_\nu\}.
$$
It can be shown that for each $\nu$,
unboundedly many members of $\bar{M}$
belongs to $B_\nu$. By maximality,
$\bar{M}\setminus B_\nu$ is
$\alpha$-finite for all $\nu<\kappa_2$.
Let $h(\nu)$ be the supremum of this
$\alpha$-finite set. We claim:
\vskip.2in
\itemitem{} For each $y\in\bar{M}$,
there are only finitely many $\nu$'s
such that $h(\nu)<y$.
\vskip.2in
Fix a $y\in\bar{M}$. Suppose there are
infinitely many $\nu$'s such that
$h(\nu)<y$. This means that $y\in
B_\nu$ for each of these $\nu$'s. Since
the $A_\nu$'s are pairwise disjoint,
$y$ must have entered the $B_\nu$'s at
different stages $\sigma$ exhibiting
infinitely many different
order types for $y\setminus M^\sigma$.
But this contradicts the
well-ordering of ordinals. This proves
the claim.

It follows from the claim that
$\kappa_2$ and hence $\alpha^*_2$ is
countable. We conclude that $\alpha$ is
countable.

\vskip.3in

\noindent{\bf  Post's  Problem Above
$\emptyset^\prime$   And  Set-theoretic Methods}
\vskip.2in
The second example in the negative
direction concerns $\alpha$-degrees
above {\bf 0$^\prime$}. We discuss how
 Silver's work on singular cardinals of
uncountable cofinality when merged with
Jensen's theory is
 exploited to derive a strong
structural difference in degree  theory
  for a class of admissible
ordinals.   Further applications are
discussed in Section II.

The problem  to consider 
is simple: Does Post's problem hold
above any $\alpha$-degree  ? In  other
words,  for any set   $A$, do there
exist sets $B$   and $C$ RE ~in $A$ such that
$A<_\alpha B$ and $A<_\alpha C$, and
$B$, $C$ have incomparable
$\alpha$-degrees ?  A  related, and
more  general,  question   asks if
there exist incomparable
$\alpha$-degrees  above  any  given degree.
A basic  theorem of Kleene-Post states
that this holds  when $\alpha=\omega$.
For $\alpha=\aleph_{\omega_1}$, the
answer turns out to be negative in a
very strong way:
\vskip.2in
\noindent
{\bf 
 Theorem 6 (Friedman [1981])} {\sl Assume V=L. If
$\alpha=\aleph_{\omega_1}$, then the
$\alpha$-degrees above
{\bf 0'} are well-ordered,
with successor provided by the jump operator.}

\noindent {\bf Proof:} A complete proof requires a heavy dose
of Jensen's fine structure theory. We
give a sketch here of the proof of
the easy half.
 Given $A, B \ge_\alpha
\emptyset^\prime$, define the growth function
$g_A$ of $A$ so that $g_A(\delta)$ is
the least ordinal $u$ such that $A\cap
\aleph_\delta\in L_u$. Define $g_B$
similarly. Then either $g_A(\delta)\ge
g_B(\delta)$ for stationarily many
$\delta$, or $g_A(\delta)<g_B(\delta)$
for closed and unboundedly many
$\delta$. Silver's analysis of growth
functions [1974]  shows that in the former
case $A\le_\alpha B$, while in the latter
case $A>_\alpha B$. As a consequence,
if $A<_\alpha B$, then
$g_A(\delta)<g_B(\delta)$ for
a closed and unbounded set of
$\delta$'s. Using this,  the well-ordering
property follows from the observation
that a countable intersection of closed
and unbounded sets is closed and
unbounded. Hence a countable descending
chain of $\alpha$-degrees above
$\emptyset^\prime$ has a least element.

In Friedman [1981]
 it is shown that the
well-ordering of $\alpha$-degrees above {\bf
0$^\prime$} is actually achieved through the
jump operator, and these $\alpha$-degrees
are represented by ``master codes'' in
Jensen's sense. 

The situation for countable cofinality
turns out to be radically different.
Harrington and Solovay have
independently shown that incomparable
$\alpha$-degrees exist above {\bf 0$^\prime$}
for $\alpha=\aleph_\omega^L$. The
following result (Chong and Mourad [in preparation])
solves Post's problem above {\bf 0$^\prime$}:
\vskip.2in
\noindent{\bf  Theorem 7} {\sl Let
$\alpha=\aleph^L_\omega$. Then there
exist sets $A$ and $B$, $\alpha$-RE in
and above
$\emptyset^\prime$, which are  of incomparable $\alpha$-degree.}
\vskip.2in
Since such sets $A$ and $B$ are
necessarily non-hyperregular, the
classical approach of finite injury
argument no longer applies. Instead, a
refinement of the method first used in
establishing the Friedberg-Muchnik
Theorem for $B\Sigma_1$-models of
arithmetic (Chong Mourad [1992]), called unions of intervals,
is exploited to ensure that all
requirements are met within
$\omega$-steps. Since $A$ and $B$ lie
above $\emptyset^\prime$ and are
therefore able to ``climb up'' 
$\alpha$ in $\omega$-many steps, the
construction succeeds.
\vskip.3in
\noindent{\bf  Applications to Fragments
of Peano Arithmetic}
\vskip.2in
One of the most interesting
applications of techniques of
$\alpha$-recursion theory in recent
years has been in
the area of {\it reverse recursion
theory\/}. Starting with the basic
axioms of Peano arithmetic without the induction
scheme, one asks: 
\vskip.2in
\itemitem{} What is the proof-theoretic
strength of a given theorem in
recursion theory ? In particular, how
much of the induction scheme is required
to prove the  theorem ?
\vskip.2in
Kirby and Paris [1978] have provided a
hierarchy of theories of increasing
proof-theoretic strength, and this
hierarchy forms the basis for the study
of subrecursive recursion theory. Let
$P^-$ be axioms of Peano arithmetic
with exponentiation but
without the induction scheme. Let
$I\Sigma_n$ denote the induction scheme
for all $\Sigma_n$ formulas, and
$B\Sigma_n$ to be replacement (collection) axiom for
$\Sigma_n$ formulas: every
$\Sigma_n$-function maps a ``finite set''
(in the sense of the given model) onto
a ``finite set''. Then with $P^-$ as the
underlying theory, one has ($n\ge 0$)
$B\Sigma_{n+1}$ to be strictly stronger than
$I\Sigma_n$, which is in turn strictly
stronger than $B\Sigma_n$.
\vskip.2in
Slaman and Woodin [1989] initiated the
study of recursion theory on fragments
of Peano arithmetic. We illustrate here
how techniques of ordinal recursion
theory are adapted to investigate
problems in this area.
\vskip.2in

\noindent{\bf Theorem 8 (Chong and Mourad
[1992])} {\sl $P^- +B\Sigma_1$ proves the
Friedberg-Mucknik theorem.}

\noindent{\bf Proof:} Simpson
(unpublished) observed that $I\Sigma_1$
was sufficient to verify that the
standard construction works. Thus let
$\Cal M$ be a model of $P^-+B\Sigma_1$
in which $\Sigma_1$-induction fails.
There is then a cofinal $\Sigma_1(\Cal
M)$ map $f$ defined on a $\Sigma_1(\Cal
M)$-definable ``cut'' $X$. This map $f$ on
$\Cal M$
acts very much like a
$\Sigma_2$-cofinal function of
$\aleph^L_\omega$ (with domain
$\omega$), or indeed a
$\Sigma_1$-cofinal map on a
rudimentarily closed $\beta$ which is
not admissible ($\beta$-recursion
theory in Section II). The idea now is
to treat  $\Cal M$ as having
``cofinality $X$'' (so that $\Cal
M=\cup_{t\in X} M_t$, and $M_t\subset M_{t+1}$), and construct a
``Friedberg-Muchnik pair''  by
satisfying the requirements
successively within each $M_x$. 
\vskip.2in
The following example shows how the
methods of Shore [1976a] is applied.
\vskip.2in
\noindent{\bf Theorem 9 (Mytilinaios
and Slaman [1988])} {\sl $P^- +
B\Sigma_2$ does not prove the existence
of an incomplete high RE~set.}

\noindent{\bf Proof:} There is a model
$\Cal M$ of $P^- +B\Sigma_2$ in which
$\Sigma_2$-induction fails (with
$\omega$ as the domain of a
$\Sigma_2(\Cal M)$-cofinal function $f$), and in
which every real is ``coded'' (meaning
it is the initial segment of a
``finite'' set). The function $f$ is
recursive in $\emptyset^\prime$. If $A$ is an
incomplete RE~set in $\Cal M$, then 
for each $n\in\omega$, there is a least
$g(n)$ such that $e\in A^\prime|f(n)$ if and only if
$\{e\}^{A^{g(n)}}_{g(n)}(e)\downarrow$,
else $A$ will be complete. Now
$n\mapsto g(n)$ is coded, and so one
may use it  to compute $A^\prime$ from $\emptyset^\prime$.
\vskip.2in
Chong and Yang
[to appear] have recently shown that the existence of a maxmal
set, as well as that of an incomplete high
set, is equivalent to $P^- +I\Sigma_2$.
In general, just as for
$\alpha$-recursion theory, infinite injury priority
method is less well understood.
Groszek, 
Mytilinaios and Slaman [to appear] have recently shown that
 $P^-+B\Sigma_2$ proves the Density
Theorem. The proof-theoretic
classification of this theorem is not known.
\vskip.2in
We refer the reader to Chong [1984] and
Sacks [1990] for more complete
treatments  on
$\alpha$-recursion theory.

\overfullrule0pt

\def\ss{\vskip4pt}
\vskip.4in

\centerline{\bf II. $\beta$-Recursion Theory}
\vskip.3in 

Studying the global structure of the $\alpha$-degrees clearly exposes the
need to deal with failures of admissibility: even though an ordinal
is admissible it may fail to be relative to a set whose degree we
wish to analyze.  Indeed, the main thrust of the work in
$\alpha$-recursion theory has been to demonstrate that
recursion-theoretic constructions from classical recursion theory
which seem to require a large amount of admissibility, say $\Sigma_2$
or even $\Sigma_3$, can actually be refined so as to succeed with only
the assumption of $\Sigma_1$-admissibility.  In view of this it is
natural to ask: Can the assumption of $\Sigma_1$-admissibility be
eliminated ?
\vskip.1in

However on hindsight it is fair to say that a stronger motivation for
the development of $\beta$-recursion theory was to find new
applications of the beautiful work of Jensen [1972] on the fine
structure of $L$, to ordinal recursion theory.  Jensen's work ignores
admissibility distinctions but concentrates only on iterations of the
jump operator (``master codes''); $\beta$-recursion theory extends
his idea to degree theory in general. 
 
\vskip.1in

The basic notions in $\beta$-recursion theory are defined using
Jensen's hierarchy for $L$, the $J_\alpha$-hierarchy, which enjoys
the following properties : 
\ss
\itemitem{(a)} $J_0 = \emptyset, J_{\alpha+1} \cap P(J_\alpha) =$
Definable subsets of $J_\alpha$ (with parameters), $J_\lambda = \cup
\{J_\alpha \mid \alpha < \lambda\}$ for limit $\lambda$. 
\ss

\itemitem{(b)} $J_\alpha$ obeys
$\Sigma_0(J_\alpha)$-comprehension and is
closed under pairing. 
\vskip.1in

Of course the improvement over the $L_\alpha$-hierarchy is closure
under pairing.  Unfortunately $J_\alpha \cap \text{ORD}$ is $\omega\alpha$
and not $\alpha$. So we define, for limit $\beta: S_\beta =
J_\alpha$ where $\beta = \omega\alpha$.  $\beta$-recursion theory
takes place on the set $S_\beta$. 
\vskip.1in

The notions $\Sigma_n$-cofinality and  $\Sigma_n$-projection apply to
$\beta$ as they do in the admissible case: $\Sigma_n$-cofinality
$(\beta) = \text{least} \gamma$ such
that there is a cofinal  $f:\gamma \to \beta$ which is 
$\Sigma_n(S_\beta)$;
$\Sigma_n$-projectum $(\beta) =
\text{least} \ \gamma$  such that
there is a one-one $f:\beta \to \gamma$ which is 
$\Sigma_n(S_\beta)$.  These are
either equal to $\beta$ or are $\beta$-cardinals (cardinals in the
sense of $S_\beta)$.  An important result of Jensen [1972] states that
$\Sigma_n$ projectum $(\beta)$ is also the least $\gamma$ such that
some $\Sigma_n(S_\beta)$ subset of
$\gamma$ is not an element of $S_\beta$. 
\vskip.1in

We are ready to define the basic notions of $\beta$-recursion theory.
As in $\alpha$-recursion theory, $A \subseteq S_\beta$ is
$\beta$-recursively enumerable,
$\beta$-recursive, $\beta$-finite if
and only if
$A$ is
$\Sigma_1(S_\beta)$, 
$\Delta_1(S_\beta)$, 
an element of $S_\beta$, respectively.  However when $\beta$ is
inadmissible (i.e., $\Sigma_1$ cofinality $(\beta) < \beta)$, a new
and stronger notion of $\beta$-RE ($\beta$-recurrsively enumerable)
arises: $A$ is {\it tamely} $\beta$-RE if 
$A^* = \{x \in S_\beta \mid x \subseteq A\}$ is $\beta$-RE.  This is
equivalent to saying that $A$ is the union of a $\beta$-recursive
sequence $\langle A^\sigma \mid \sigma < \beta \rangle$ with the
property that if $x \subseteq A$ and  $\beta$-finite then $x \subseteq
A^\sigma$ for some $\sigma < \beta$. 
\vskip.1in

The weak and strong reducibilities $\leq_{w\beta}$, $\leq_{\beta}$ are
defined as they are in $\alpha$-recursion theory: One way of
achieving these definitions is through the use of ``neighborhood
conditions'': define $N(A) = \{\langle x,y \rangle \mid x,y \;\hbox{are}\;
\beta-$finite, $x \subseteq A, y
\subseteq \bar{A}\}$.   $B$ is $\beta$-RE in $A$ if for some
$\beta$-RE $W$, $x\in B$ if and only if
$\exists z \in N(A)\;
[(x,z) \in W]$. Then $B \leq_{w\beta}A$
if and only if $B$,
$\bar{B}$ are both $\beta$-RE in $A$,
and $B \leq_{\beta} A$ if and only if 
$B^*$ and $\bar{B}^*$ are both
$\beta$-RE in $A$ (if and only if $B$, $\bar{B}$
are both ``tamely'' $\beta$-RE in $A$).  The strong reducibility
$\leq_\beta$ is transitive. 
\vskip.1in

Now some genuinely new phenomena arise in the inadmissible case, with
regard to $\beta$-reducibility.  These are summarized in the 
following result. 
\vskip.3in

\noindent {\bf Theorem 1 (Friedman [1979])}
{\sl Assume that $\beta$ is
inadmissible.   Then there is a $\beta$-recursive set $A$ such that: 
\itemitem{(i)} $\emptyset <_\beta A
<_\beta C$ where $C$ is a  complete $\beta$-RE
set. 
\itemitem{(ii)} Any tamely $\beta$-RE set and any $\beta$-recursive set
is $\beta$-reducible to $A$. 
\itemitem{(iii)} $C \leq_{w\beta} A.$}

\vskip.1in

\noindent Thus $\beta$-recursiveness does not imply $\beta$-reducibility to
$\emptyset$, and the complete $\beta$-RE set is weakly
$\beta$-reducible to a $\beta$-recursive set !
\vskip.1in

It is easy to define $A$ (in fact $A$ can be taken to be a
$\Delta_1$ master code in the sense of Jensen [1972]). 
Let $f:\Sigma_1-\text{cofinality} (\beta) \to \beta$ be 
$\Sigma_1(S_\beta)$ and cofinal,
and take $A = \{(e,x,\gamma)\mid x \in W_e$ by stage 
$f(\gamma), \gamma < \Sigma_1$-cofinality $(\beta)\}$, where $W_e$
is the $e$th $\beta$-RE set.  Then $A$ is $\beta$-recursive and since
$x \not\in W_e$ if and only if $\{e\}\times \{x\}\times \gamma \subseteq
\bar{A}$ we get $C \leq_{w\beta} A$.  The other properties are
not difficult to verify. 
\vskip.1in

The $\beta$-degree of $A$ is referred
to as {\bf 0}$^{1/2}$ and serves as a
new type of jump operator in $\beta$-recursion theory.  Of course
{\bf 0}$^{1/2}$ provides an easy solution to a version of Post's Problem in
the inadmissible case; however it does not answer the following
question, which has come to be adopted as the official version of
Post's Problem in $\beta$-recursion theory. 
\vskip.1in

\noindent {\bf Post's Problem} \ \ Do there exist $\beta$-RE sets $A, B$
such that $ A \not\leq_{w\beta} B, \; B \not\leq_{w\beta} A$ ? 
\vskip.1in

As in $\alpha$-recursion theory, Post's Problem has served as a
driving force behind much of the work in $\beta$-recursion theory. 
\vskip.1in

Early on it became apparent that with regard to questions such as
Post's Problem the inadmissible ordinals divide into two very
different classes.  (This distinction occurred earlier in Jensen's
proof of $\Sigma_2$ uniformization.)
$\beta$ is {\it weakly
admissible} if $\Sigma_1$-cofinality
$(\beta)  \geq \Sigma_1$-projectum
$(\beta)$.  Otherwise $\beta$ is {\it strongly
inadmissible}.  In the former case many arguments from
$\alpha$-recursion theory can be adapted, for the following reason :
if $\beta$ is weakly admissible (but inadmissible) then there is a
$\beta$-recursive bijection between
$S_\beta$ and  $\Sigma_1$-cofinality
$(\beta)$.  Moreover there  is a $\beta$-recursive $A
\subseteq \Sigma_1$-cofinality $(\beta) = \kappa$ which is a $\Delta_1$
master code for $S_\beta$ in Jensen's sense: $B \subseteq K$ is
$\beta$-RE iff $B$ is 
$\Sigma_1(L_\kappa, A)$.  Thus $\beta$-recursion theory is
closely related to $\kappa$-recursion theory, relativized to $A$ and the
structure $(L_\kappa, A)$ is admissible.  This is sufficient
to reduce the solution to Post's Problem for $\beta$ to the
previosuly known (positive) solution for 
$(L_\kappa, A)$.  Exactly how much can be reduced from
$\beta$ to $(L_\kappa, A)$ is analyzed in 
Maass [1978a]. 
\vskip.1in

The greater challenges in $\beta$-recursion theory arise in the
strongly inadmissible case.  Techniques from admissibility theory no
longer apply; instead methods from combinatorial set theory are
needed. The first attack on Post's Problem in the strongly
inadmissible case appears in Friedman [1980]. 
\vskip.1in

\noindent{\bf Theorem 2 (Friedman [1980])} Suppose $\beta$ has regular
projectum: $\Sigma_1$-projectum $(\beta)$ is regular with respect to
$\beta$-recursive functions.  Then Post's Problem has a positive
solution. 
\vskip.1in

The proof uses an effective analog of Jensen's
$\diamondsuit$-principle.  We provide here a sketch of the proof, in
the special case where $\Sigma_1$-projectum $(\beta) = \aleph_1^L$.
We may assume that $\Sigma_1$-confinality $(\beta) = \omega$ (else
$\beta$ is weakly admissible) but actually the proof makes no use of
this. 
\vskip.1in

We build $\beta$-RE $A, B \subseteq \aleph_1^L$ so as to meet the
requirements $R_e^A : \bar{B} \not= W_e^A$ and 
$R_e^B : \bar{A} \not= W_e^B$, where of course $W_e^A$ is the
$e$th set $\beta$-RE in $A$.  To achieve $R_e^B$ we want an $x
\not\in A$ and a neighborhood condition $y \subseteq B, \; z
\subseteq \bar{B}$ so that $(x,(y,t))\in W_e$.  One difference from the admissible case is that we
may in fact have to actively guarantee $y \subseteq B$ as otherwise
there may be no stage $\sigma < \beta$ where $y \subseteq B^\sigma$,
due to the lack of tameness.  It is possible however to arrange a
weak form of tameness (through use of additional requirements) to
insure that in fact $y - B^\sigma$ is countable at some stage, so we
need only act to put a countable set into $A$ or $B$ for the sake of
each requirement. 
\vskip.1in

The second and most striking difference from the admissible case is
that we act on each requirement {\it at most once}.  What
enables us to make this restriction is the following.  Requirements
can be listed in a sequence $\langle R_\delta\mid\delta <
\aleph_1^L \rangle$ and as we are only putting countable sets into
$A$ or $B$ to satisfy requirements there will be a closed unbounded
set of requirements $R_\delta$ such that all action  taken by
$R_\delta'$, $\delta' < \delta$ takes place below $\delta$.
Moreover $R_\delta$ will only seek to protect ordinals $\geq \delta$
from entering $A$ or $B$ so will never be injured.  If we arrange
that each requirement appears as $R_\delta$ for a stationary set of
$\delta$'s then each requirement will have the opportunity to act
without injury.  (So in fact this not really an injury argument at
all.) 
\vskip.1in

Finally notice however that we have prohibited requirement $R_\delta$
from taking any action below $\delta$; this  requires that $R_\delta$
has a way of ``guessing'' at $A \cap \delta, \; B \cap \delta$.  The
necessary guesses are provided by Jensen's $\diamondsuit$-principle.
We end this sketch with no more than a statement of $\diamondsuit$.
\vskip.2in

\noindent ({\bf $\diamondsuit$-principle}) Suppose $E \subseteq \aleph_1^L$ is
stationary. Then there exists $\langle G_\alpha\mid\alpha \in E
\rangle$ such that : 
\ss
\noindent{(a)} $G_\alpha \subseteq \alpha$ for $\alpha \in E$. 

\item {(b)} If $A \subseteq \aleph_1^L$ then $\{\alpha \in E \mid A
\cap \alpha = G_\alpha\}$ is stationary. 
\ss

\noindent (In the general case of Theorem 2 we must weaken this somewhat
but the general idea is the same.)
\vskip.1in

The final case, where $\beta$ is strongly inadmissible with
{\it singular} projectum is entirely different.  In fact Post's
Problem may have a negative solution!  We illustrate the result with
a typical example : $\beta = \alpha\cdot\omega$ where $\alpha =
\aleph_{w_1}^L$. 
\vskip.1in

\noindent {\bf Theorem 3 (Friedman
[1978])} {\sl Let $C$ be the complete
$\beta$-RE set.  If $A$ is $\beta$-RE then either $A \leq_\beta
\emptyset$ or $C \leq_{w\beta} A$.} 
\vskip.1in

The proof makes use of the work in Silver [1974] on the singular
cardinal problem in set theory (as was
for Theorem 5.1 in Section I).  We
confine ourselves here to only a
very rough sketch of the proof.  The main idea is to look at
{\it growth rates} for subsets of $\alpha$.  Specifically,
suppose $A \subseteq \alpha$ is constructible and define $f_A(\alpha)
= \hbox{least} \; \delta$ such that $A \cap \aleph_\gamma^L$ belongs
to $L_\delta$. Then it can be shown
that if $f_A(\gamma) \leq f_B (\gamma)$ for unboundedly many $\gamma$
then in fact $A$ is weakly $\beta$-reducible to $B$.  This can be
extended to $\beta = \alpha\cdot\omega$ to show that in fact any two
subsets of $\beta$ are $\leq_{w\beta}$-comparable.  If $C$ is the
complete $\beta$-RE set then associated to $C$ is a growth rate $f$
which is the limit of $\beta$-finite growth rates $f_n, \; n\in\omega$.
Thus either $f_A$ is dominated by some $f_n$ and is hence
$\beta$-finite or $f_A$ dominates $f$ in which case $C \leq_{w\beta}
A$.  The uncountable cofinality of $\alpha$ is used both to apply
Silver's work and to simultaneously bound the $f_n$'s in this last
argument. 
\vskip.1in

\vskip.4truein

\centerline {\bf III. The Admissibility Spectrum}

\vskip.3in

Until now we have fixed an ordinal $\alpha$ (admissible or not) and
studied definability for subsets of $\alpha$.  In this section we
invert the process: fix a subset $x$ of some cardinal $\kappa$, a
theory $T$ and consider the $T$-{\it spectrum} of $x$ =
$\Lambda_T(x) = \{\alpha \mid L_\alpha[x] \models T\}$.  Thus natural
classes of ordinals can be defined from sets $x$ and we can ask for a
characterization of which classes arise in this way. 
\vskip.1in

Most of the work in this area has concentrated on the case $\kappa =
\omega, \; T = KP $= Admissible Set Theory.  However there is a good
understanding of $\alpha_T(x) $ = min $\Lambda_T (x)$ for arbitrary
$\kappa$ and other theories such as $KP_n = \Sigma_n$-Admissibility,
$ZF$. We will mention some of the latter work as well. 
\vskip.1in

The first result in this area is due to Sacks. 
\vskip.1in

\noindent {\bf Theorem 1 (Sacks [1976])} {\sl If $\alpha > \omega$ is
admissible and countable then $\alpha = \omega_1^R = \alpha_{KP}(R)$
for some real $R$. }
\vskip.1in

There are many proofs of 
Theorem 1, but the most adaptable (see
Friedman [1986]) is via the
method of almost disjoint forcing.  As a first
step we can add $A_0 \subseteq \alpha$ so that $\alpha$ is
$A_0$-admissible and $L_\alpha [A_0] \models $ every set is
countable. This is done by (Levy) forcing with finite conditions
${\cal P}$  from $\alpha \times \omega$ into $\alpha$ such that
${\cal P}(\beta,n) < \beta$.  Second, we can add $A_1 \subseteq
\alpha$ so that $\alpha$ is $(A_0, A_1)$-admissible and
$\beta < \alpha$ implies $\beta$ is not $(A_0 \cap \beta, A_1 \cap
\beta)$-admissible.  This is done with
conditions ${\cal P}:\beta \to 2$  such
that $\beta' \leq \beta$ implies $\beta'$ is not 
$(A_0 \cap \beta', {\cal P} \cap \beta')$-admissible. 
\vskip.1in

Now we can canonically assign a real $R_\beta$ to each $\beta <
\alpha$ so that if $\beta_1 \not=
\beta_2$ then $R_{\beta_1} \cap
R_{\beta_2}$ is finite.  By ``canonical'' we mean that $R_\beta$ is
defined in $L_{\beta+1}[A_0 \cap \beta, A_1 \cap \beta]$, uniformly.
Then we code $A = A_0 \vee A_1$ by a real $R$ using conditions
$(r,\overline{r})$ where $r$ is a finite subset of $\omega,
\overline{r}$ a finite subset of $\{R_\beta \mid \beta \in A\}$ and 
$(r_0, \overline{r}_0) \leq (r_1, \overline{r}_1)$ if 
$r_0 \supseteq r_1, \; \overline{r}_0 \supseteq \overline{r}_1$ and $n \in r_0 -
r_1$ implies  $n \not\in R_\beta$ for each $R_\beta \in \overline{r}_1$.
The result is that $\beta \in A$ if and
only if $R \cap R_\beta$ is finite and
thus $A \cap B$ is $\Delta_1$ over $L_\beta [R]$ for each $\beta <
\alpha$.  So $\beta$ is
not $R$-admissible for $\beta < \alpha$.
Preserving the admissibility of
$\alpha$ requires a bit of care, but
is based on the simple fact that almost disjoint forcing satisfies
the countable chain condition. 
\vskip.1in

Jensen extended Sacks' result to countable sequences of countable
admissibles.  For a proof of the following result see Friedman [1986].
\vskip.1in

\noindent {\bf Theorem 2 (Jensen)} {\sl Supose $X$ is a countable set of
countable admissibles greater than $\omega$ and $\alpha \in X \to
\alpha$ is $X \cap \alpha$-admissible.  Then for some real $R$, $X$
is an initial segment of $\Lambda_{KP}(R)$.} 
\vskip.1in

The proof strategy for Theorem 2 is similar to that used in Theorem
1: first add $A \subseteq \alpha$ preserving admissibility so that
$\beta < \alpha$ is $A \cap
\beta$-admissible if and only if $\beta \in X$, and
then code $A$ by a real using almost disjoint forcing.  But as we
must preserve the admissibility of ordinals in $X$ (while destroying
admissibility for ordinals not in $X$) the argument is more delicate
and has the interesting feature that extendibility of conditions 
for
the desired forcing is established using forcing. 
\vskip.1in

There are severe limitations on how much more can be done concerning
admissibility spectra in $ZFC$ alone.  This is illustrated by the
next result.  A class $X \subseteq \text{ORD}$ is $\Sigma_1$-complete if $Y$
is $\Delta_1([X], X)$ whenever $Y \subseteq \text{ORD}$ is
$\Sigma_1(L)$. 
\vskip.1in

\noindent {\bf Theorem 3} {\sl Let $\Lambda (R)$ abbreviate
$\Lambda_{KP}(R) = \{\alpha\mid\alpha \hbox{~is~$R$-admissible}\}$.
\ss
\noindent{(a)} $R \in L \to \Lambda(R) \supseteq \Lambda(0) - \beta$ for
some $\beta < \aleph_1^L$. 

\noindent{(b)} If $R \in L[G]$, $G$ is ${\cal
P}$-generic over $L$, and ${\cal P} \in
L$, then
$\Lambda (R) 
\supseteq \Lambda(0) - \beta$ for some $\beta$. 

\noindent{(c)} Suppose that $R \in
L[G]$ and $G
\subseteq {\cal P}$. If $G$ is ${\cal P}$-generic
over the amenable structure $(L,
{\cal P})$, then $\Lambda
(R)$ is not $\Sigma_1$-complete. }
\vskip.1in

\noindent {\bf Proof:} \quad (a) Let $\beta$ be large enough so that $R \in
L_\beta$. (b) Let $\beta$ be large enough so that ${\cal P} \in
L_\beta$. (c) If $\Lambda (R)$ is $\Sigma_1$-complete then $L\text{-Card} =
\{\kappa\mid L \models \kappa \hbox{~is~a~cardinal~}\}$ is 
$\Delta_1(L[R])$ and hence by
reflection, $(\kappa^+)^L < \kappa^+$ for large enough cardinals
$\kappa$.  By Jensen's Covering Theorem, $0^{\#} \in L[R]$.  But
$0^{\#}$ does not satisfy the hypothesis of (c) (see
Beller, Jensen and Welch [1982]). 
\vskip.1in

By (a), (b) of this result we see that class-forcing is required to
get a nontrivial admissibility spectrum (without assuming $0^{\#})$ 
and we should not expect such a spectrum to be $\Sigma_1$-complete. 
\vskip.1in

Using a variant of Jensen coding, R. David and S. Friedman
independently obtained a class-generic real $R$ such that
$\Lambda_{KP}(R) \subseteq$ Admissible Limits of Admissibles.  This
is a special case of the following result which appeared in David [1989].
\vskip.1in

\noindent {\bf Theorem 4 (David and Friedman)}
{\sl Suppose $\varphi (\alpha)$
is $\Sigma_1$ and $\alpha \in L\hbox{-Card} \to L \models \varphi
(\alpha)$.  Then there is a real $R$ class-generic over $L$ such that
$\Lambda_{KP}(R) \subseteq \{\alpha \mid L \models \varphi
(\alpha)\}$. }
\vskip.1in

This result is optimal in the sense that if $\varphi(\alpha)$ is the
$\Pi_1$ formula ``$\alpha$ is a cardinal'' then the conclusion must
fail by Theorem 3(c). 
\vskip.1in

We give some idea of the proof of Theorem 4.  The desired forcing is
made up of certain ``building blocks'' that are not difficult to
describe. Jensen coding methods are used to put these building blocks
together. 
\vskip.1in

We wish to arrange that if $\alpha$ is
$R$-admissible then  $\alpha$ is a limit
of admissibles.  Suppose that we have
$D \subseteq \aleph_1^L$ so that
if  $\alpha$ is $D$-admissible then
$\alpha$ is a limit of admissibles.  Then
we could hope to choose $R$ so as to code $D$ and satisfy the desired
property. 
\vskip.1in

The problem is that if we code $D$ by $R$ in the usual way (with
almost disjoint forcing) we only get: for all $\alpha, \ D \cap
(\aleph_1)^{L_\alpha}$ is $\Delta_1(L_\alpha [R])$.  So in fact
what we need about $D$ is: $L_\alpha [D \cap \xi] \models KP + \xi =
\aleph_1$ implies $\alpha$ is a limit of admissibles.  For then we need
only recover $D \cap (\aleph_1)^{L_\alpha}$ inside $L_\alpha [R]$ to
guarantee that $\alpha$ is a limit of admissibles. 
\vskip.1in

How do we obtain $D$ ? The natural thing is to force with conditions
$d$ which are initial segments of a potential $D$.  Now we come to
the main points in the proof. 
\ss
\noindent{(1)} Extendibility is easy for this forcing because given $d$
and $\gamma < \aleph_1^L$ we are free to extend $d$ to length
$\gamma$ by killing the admissibility of all ordinals between
$\sup(d)$ and $\gamma$.  It is crucial for this argument that we are
only concerned with killing admissibility, not with preserving it. 

\noindent{(2)} Cardinal-preservation for this forcing is easy to prove
assuming there is $D_2 \subseteq (\aleph_2)^L$ such that: $L_\alpha [D
\cap \xi] \models KP + \xi = \aleph_2 \to \alpha$ a limit of
admissibles. 
\vskip.1in

Thus we are faced with the original
problem, one cardinal higher !  The
solution (due to Jensen in the proof of his Coding Theorem) is to
build $R, D_1, D_2, \cdots $ simultaneously. 
\vskip.1in

Finally we introduce the requirement of  admissibility preservation
into the above.  Note that in the conclusion of Theorem 4 we have 
$\subseteq$ and not equality; indeed the freedom to kill
admissibility is crucial to the extendibility argument in (1) above. 
\vskip.1in

Nonetheless we can ask for a real $R$ for which we can control its
(nontrivial) admissibility spectrum.  This requires the method of
strong coding. 
\vskip.1in

\noindent {\bf Theorem 5} (Friedman [1987]) There is a real $R$,
class-generic over $L$, such that $\Lambda_{KP}(R) =$ Admissible
limits of admissibles. 
\vskip.1in

To prove Theorem 5 we can approach the problem much as in the proof
of Theorem 4, however extendibility of conditions is much more
difficult.  The desired extension of $d$ to length $\gamma$ must be
made generically, so as to preserve the admissibility of admissible
limits of admissibles.  (Note that this idea was foreshadowed by
Jensen's proof of Theorem 2.)  Thus conditions must be constructed
out of generic sets for ``local'' versions of the very same forcing.
So in fact we construct a strong coding ${\cal P}^\beta \subseteq
L_\beta$ at each admissible $\beta$ and then inductively build ${\cal
P}^\beta$ out of generic sets for various ${\cal P}^{\beta'}, \
\beta' < \beta$. 
\vskip.1in

A complete characterization of admissibility spectra is not known.  A
related question, which may indeed be a prerequisite for such a
characterization, is the following : which $A \subseteq \hbox{ORD}$ can be
$\Delta_1$-definable in a real class-generic over $L$ ?  On this
latter problem there has been some
significant progress. The following
will appear in Friedman-Velickovic [1995]. 
\vskip.1in

\noindent {\bf Theorem 6 (Friedman)} Suppose $V = L$
and that $A \subseteq \text{ORD}$ obeys the Condensation Condition.  Then
$A$ is $\Delta_1$ in a real class-generic over $L$, preserving
cardinals. 
\vskip.1in

We refer the reader to Friedman-Velickovic [1995] for a definition of
the Condensation Condition and a proof of Theorem 6. 
\vskip.2in

\noindent {\bf Other Work}
Much is known about $\alpha_T(x) = \min \Lambda_T(x), \ x \subseteq
\kappa$, for $T = KP_n, \ ZF$ and arbitrary infinite cardinals
$\kappa$, assuming $V = L$.  We confine ourselves here to only a
brief account. 
\vskip.1in

First we consider the (remaining) cases when $\kappa = \omega$. 
\vskip.1in

\noindent {\bf Theorem 7} 
{\sl \noindent{(a)} {\rm (Sacks [1976]})
$\alpha_{KP_n}(R)$, $R\subseteq\omega$,
can be any countable  $\Sigma_n$-admissible ordinal
greater than $\omega$. 

\noindent{(b)} {\rm (David [1982],
Beller in Beller, Jensen and Welch [1982])}
$\alpha_{ZF}(R)$, $R\subseteq\omega$, can be any countable $\alpha$
such that $L_\alpha \models ZF$. }
\vskip.1in

Theorem 7(a) can be proved much like Theorem 1.  For Theorem 7(b)
note that it suffices to first find
$R_0$ such that $\beta$ an
$L$-cardinal implies that $L_\beta[R_0]\not\models
ZF$ and then apply Theorem 4 (relativized to $R_0$).  The former is
not hard to arrange using only the statement of Jensen's Coding
Theorem. (Of course historically Theorem 7(b) was proved directly as
Theorem 4 was not available.) 
\vskip.1in

Next suppose that $\kappa$ is regular and uncountable. 
\vskip.1in

\noindent {\bf Theorem 8 (Friedman [1982])} {\sl $\alpha=\alpha_{KP_n}(x)$
for some $x\subseteq\kappa$ if and
only if $\alpha<\kappa^+, \ \alpha$ is
$\Sigma_n$-admissible, cofinality $(\alpha)=\kappa$ and $L_\alpha$ is
closed under the function $\beta\mapsto\beta^{<\kappa}$.}
\vskip.1in

The difficult part of Theorem 8 is the necessity of the stated
condition, which draws heavily on Jensen's fine structure theory.
The sufficiency is based on an almost disjoint forcing argument, not
unlike Theorem 7(a). 
\vskip.1in

\noindent {\bf Theorem 9 (David and Friedman [1985])} {\sl
 $\alpha=\alpha_{ZF}(x)$ for some $x\subseteq\kappa$ if
and only if $\kappa<\alpha<\kappa^+$,
$L_\alpha\models ZF$  and there are
$\beta<\alpha$,  $\langle X_n\mid n\in\omega\rangle$ such that 
\itemitem{(i)} $\forall\gamma<\kappa\forall
f:\gamma\to\beta$  ($f$ bounded $\to
f\in L_\alpha)$,
\itemitem{(ii)} $X_n \in
L_\alpha$, $L_\alpha$-Card $(X_n)$ is
less than $\beta$ for all $n$, and $L_\alpha = \cup
\{X_n \mid n\in\omega\}$, and
\itemitem{(iii)} $\beta$ is a regular cardinal in
$L_\alpha$.} 
\vskip.1in

The proof of Theorem 9 makes use of almost disjoint forcing, the
Covering Theorem (relativized to some $L_\alpha [x])$ and Jensen's
fine structure theory. 
\vskip.1in

When $\kappa$ is singular of cofinality $\omega$ then methods from
infinitary model theory come into play. 
\vskip.1in

\noindent {\bf Theorem 10 (Friedman [1981a])}
{\sl $\alpha=\alpha_{KP}(x)$
for some $x\subseteq\kappa$ if and only
if 
\itemitem{(i)} $\kappa<\alpha<\kappa^+$,
\itemitem{(ii)} if
there is a largest $L_\alpha$-cardinal $\gamma$ then cofinality
$(\gamma)=\omega$, and
\itemitem{(iii)} there is a 1-1 function $f: L_\alpha\to
\kappa$ such that $f^{-1}[\delta]\in L_\alpha$ for each $\delta<
\kappa$. }
\vskip.1in

Under the conditions stated in Theorem 10, a version of the Barwise
Compactness Theorem is established, which can then be used to obtain
the desired $x$.  A related result
appears in Magidor, Shelah and 
Stavi
[1984]. 
\vskip.1in

For $n>1$ a surprising thing occurs: for any $x\subseteq
\aleph_\omega$, $x\in
L_{\alpha_{KP_2}}$ ! And an even stronger fact
holds for $x\subseteq\aleph_{\omega_1}$, namely $x\in
L_{\alpha_{KP}}(x)$. Both of these facts follow from an effective
version of Jensen's Covering Theorem.  This puts severe restrictions
on the possible values for
$\alpha_{KP_n}(x)$, $x\subseteq \kappa$ for
$n>1$, $\kappa$ singular of cofinality $\omega$ and for $n\geq 1$, $\kappa$ singular
of uncountable cofinality (as well as for $\alpha_{ZF}(x))$.  The
reader is referred to Friedman [1981],
David and Friedman [1985] for complete
characterizations. 
\vskip.1in

\vskip.3in

\vskip.1in

\parindent=3truecm
\centerline{\bf References}
\vskip.2in
\noindent{A.~Beller, R.~B.~Jensen, and P.~Welch [1982]}   {\it Coding the Universe\/},
London Mathematical Society Lecture
Notes Vol. 47

\noindent C.~T.~Chong [1984], {\it Techniques
of Admissible Recursion Theory\/},
Lecture Notes in Math. Vol 1106,
Springer Verlag

\noindent{C.~T.~Chong and K.~J.~Mourad} [1992],
$\Sigma_n$-definability without
$\Sigma_n$-induction, {\it
Trans.~AMS\/} Vol.~334, 349--363

\noindent{C.~T.~Chong and K.~J.~Mourad} [in preparation], Post's
problem and singularity

\noindent{C.~T.~Chong and Y.~Yang} [to
appear], Maximal sets, high RE sets,
and $\Sigma_2$-induction

\noindent{R.~David [1982]}   Some applications of Jensen's Coding Theorem,
{\it Ann.~Math.~Logic\/} Vol.~22,  177-196 

\noindent{R.~David and S.~D.~Friedman [1985]}    Uncountable
$ZF$-ordinals, {\it Proc. Symposia in
Pure Math.\/}  AMS, Vol.~42, 217-222.

\noindent{R.~David [1989]}   A functional
$\Pi_2^1$-singleton, {\it Advances in
Math.\/} Vol.~74,  258-268

\noindent{S.~D.~Friedman [1978]}  Negative solutions to Post's Problem, I, in
{\it Generalized Recursion Theory II},
Fenstad-Gandy-Sacks,  Eds., North-Holland 

\noindent{S.~D.~Friedman [1979]}
$\beta$-recursion theory, {\it Transactions
AMS\/} Vol.~255, 173-200

\noindent{S.~D.~Friedman} [1980]  Post's Problem without admissibility,
{\it Advances in Math.\/} Vol.~35, 30-49

\noindent{S.~D.~Friedman} [1981], Negative
solutions to Post's problem, II. {\it
Ann.~Math.\/} Vol.~113, 25--43

\noindent{S.~D.~Friedman [1981a]}   Uncountable admissibles, II: Compactness,
{\it Israel  J.~Math.\/} Vol.~40, 129-149

\noindent{S.~D.~Friedman [1982]}  Uncountable admissibles, I: Forcing,
{\it Trans.~AMS\/} Vol.~270, 61-73

\noindent{S.~D.~Friedman [1983]}   Some recent developments in higher recursion
theory, {\it J.~Symbolic Logic\/} Vol.~48, 629-642

\noindent{S.~D.~Friedman [1986]}   An introduction to the admissibility spectrum,
{\it Logic, Metholology and Philosophy
of Science VII\/},
Marcus-Dorn-Weingartner, Eds., North-Holland,  129-139 

\noindent{S.~D.~Friedman [1987]} Strong coding, {\it Ann. Pure and Applied Logic\/}
Vol. , 1-98 

\noindent{S.~D.~Friedman and Velickovic [1995]}
$\Delta_1$-definability, {\it to
appear\/} 

\noindent{M.~Groszek, 
M.~Mytilinaios and T.~A.~Slaman}
[to appear], The Sacks density theorem
and $\Sigma_2$ bounding

\noindent{R.~B.~Jensen} [1972], The fine structure of the
constructible universe, {\it
Ann.~Math.~Logic\/} Vol.~4, 229--308

\noindent{S.~C.~Kleene} [1938], On notation for
ordinal numbers, {\it J.~Symbolic Logic\/} Vol.~3, 150--155

\noindent{G.~Kreisel and G.~E.~Sacks} [1965],
Metarecursive sets, {\it J.~Symbolic
Logic\/} Vol.~10, 318--336

\noindent{M.~Lerman} [1972], On suborderings of
the $\alpha$-recursively enumerable
degrees, {\it Ann.~Math. Logic\/}
Vol.~4, 369--392

\noindent{M.~Lerman} [1974], Maximal
$\alpha$-RE~sets, {\it Trans.~AMS\/}
Vol.~188, 341-386

\noindent{W.~Maass} [1978], The uniform regulary set
theorem in $\alpha$-recursion theory,
{\it J.~Symbolic Logic\/} Vol.~43, 270--279

\noindent{W.~A.~Maass [1978a]}
Inadmissibility, tame RE sets and the admissible
collapse, {\it Ann.~Math. Logic\/} Vol.~13, 149-170

\noindent{M.~Magidor, S.~Shelah, and
J.~Stavi [1984]}    Countably decomposable admissible
sets, {\it Ann. Math.~Logic\/} Vol.~26,
287--361  

\noindent{M.~Mytilinaios and T.~A.~Slaman}
[1988], $\Sigma_2$-collection and the
infinite injury priority method, {\it
J.~Symbolic Logic\/} Vol.~53, 212--221

\noindent{J.~B.~Paris and L.~A.~Kirby}
[1978], $\Sigma_n$ collection schemas
in models of arithmetic, in: {\it Logic
Colloquium '77\/}, North-Holland

\noindent{H.~Rogers Jr.} [1967], {\it Theory of
Recursive Functions and Effective
Computability\/}, McGraw-Hill

\noindent{G.~E.~Sacks} [1966], Post's
problem, admissible ordinals, and
regularity, {\it Trans.~AMS\/}
Vol.~124, 1--23

\noindent{G.~E.~Sacks [1976]}   Countable admissible ordinals and hyperdegrees,
{\it Advances in Math.\/} Vol.~20, 231-262 

\noindent G.~E.~Sacks [1990], {\it Higher
Recursion Theory\/}, Springer Verlag

\noindent{G.~E.~Sacks and S.~G.~Simpson} [1972],
The $\alpha$-finite injury method, {\it
Ann.~Math.~Logic\/} Vol.~4, 323--367 

\noindent{R.~A.~Shore} [1976], The recursively enumerable
$\alpha$-degrees are dense, {\it
Ann.~Math.~Logic\/} Vol.~9, 123--155

\noindent{R.~A.~Shore} [1976a], On the jump of an
$\alpha$-recursively enumerable set,
{\it Trans.~AMS.\/} Vol.~217, 351--363 

\noindent{J.~H.~Silver} [1974], On the singular
cardinals problem, {\it
Proc.~International Congress of
Mathematicians 1974\/}, 265--268

\noindent{T.~A.~Slaman and W.~H.~Woodin} 
[1989], $\Sigma_1$-collection and the
finite injury priority method, in: {\it
Mathematical Logic and Its
Applications\/}, Lecture Notes in
Mathematics Vol.~1388, Springer Verlag

\bye